\documentclass[letterpaper]{article}%
\usepackage{graphicx}
\usepackage{amssymb}
\usepackage{amsmath}
\usepackage{amsfonts}
\usepackage{sw20bams}
\usepackage{comment}
\usepackage{hyperref}
\usepackage{url}
\usepackage{aurl}%
\providecommand{\U}[1]{\protect\rule{.1in}{.1in}}

\ifx\pdfoutput\relax\let\pdfoutput=\undefined\fi
\newcount\msipdfoutput
\ifx\pdfoutput\undefined\else
\ifcase\pdfoutput\else
\ifx\paperwidth\undefined\else
\ifdim\paperheight=0pt\relax\else\pdfpageheight\paperheight\fi
\ifdim\paperwidth=0pt\relax\else\pdfpagewidth\paperwidth\fi
\fi\fi\fi
\begin{document}

\title{Fast Ramanujan--type Series for Logarithms. Part I.}
\author{Jorge Zuniga\thanks{$\ \ $Independent Researcher. \texttt{jorge.zuniga.hansen.at.gmail.com}}}
\date{May 2025\footnote{$\,\ $Keywords: Number Theory, Sequences \& Series,
Algorithms, Binary Splitting, Wilf-Zeilberger, LLL}\vspace*{-4pt}}
\maketitle

\begin{abstract}
This report introduces new series and variations of some hypergeometric type
identities for fast computing of logarithms log$\,p$ for small positive
integers $p.$ These series were found using Wilf--Zeilberger (WZ) and/or
integer detection algorithms (LLL) providing highly efficient linearly
convergent rational approximants for these constants. Some of the new
identities are of {}$_{4}F_{3}$ type, but higher ones are found as well and
hypergeometric series for log$\,p$ (with variable $p$) have been derived.
Found identities are proven by \textit{i) }classical Beta Integral methods,
\textit{ii)} some hypergometric closed forms and\textit{\ iii)} rational
certificates from the WZ method. Since they are very fast, these series are
particularly suitable to be embodied in mathematical software being
implemented as binary splitting which produces very efficient algorithms. Over
$10^{12}$\thinspace decimal places have been obtained for some logarithms in
reasonable time.

\end{abstract}

\section{Introduction\smallskip}

Traditionally the calculation of logarithms\ of natural numbers, including
small primes, has been performed by means of highly efficient hypergeometric
Machin (MT) formulas that are linear combinations of some fast convergent
arcotanh series (\cite{FRED}\thinspace Table 1). In this work, by applying
diverse searching strategies, many members of a new infinite family of
logarithm identities based on a different kind of hypergeometric series were
found. From this family, only the identities that are more efficient than MT
formulas were selected. For efficiency, it is understood the lesser computing
time per digit that a particular hypergeometric series delivers if the binary
splitting algorithm is applied. See \cite{YAK} and references therein. An
exact definition goes through Eq.(\ref{3}), the concept of \textit{Binary
Splitting Cost}.

\bigskip In this searching process the LLL lattice reduction algorithm
\cite{LLL} was intensively used to look for integer relationships between the
logarithm constant\ probe and certain numerically evaluated partial series. In
this case the selected relationships found are just conjectured identities
that require independent proofs. Fortunatedly most of these new discovered
identities are independently proven and generalized by using\ Beta
function\ techniques (Gourevitch--Guillera, \cite{GUIGOU}, Section 5) and also
by applying some proofs of hypergeometric sums recently found by Z.W. Sun and
Y. Zhou \cite{SUNZHOU} and C. Hakimoglu \cite{CHBROWN}. This last work was
applied to derive and generate the whole family of these\ level 1 logarithm
series Eq.(\ref{18z}), where some of its members were also found by the
searching process.

\bigskip On the other hand, the Wilf--Zeilberger (WZ) method \cite{WILFZEIL}%
\cite{KOEPF} was also applied to obtain new proven formulas (of Ramanujan-type
and beyond) where some of them were also found by the LLL algorithm. In these
cases a generalization of the technique used by Amdeberhan and Zeilberger
\cite{ZEIL} was applied allowing that all the selected series get their
automated proofs directly from the algorithm since it also delivers the
corresponding rational certificates of the new discovered identities.

\section{Hypergeometric-Type Sums\medskip}\label{sec:hypergeometric-type-sumsvspace4pt}

Let $p(n)$, $q(n)$ and $r(n)$ be integer coefficient polynomials such that
none of them vanishes for$\ n\in%
\mathbb{N}
$, and polynomials $r(n)$ and $q(n)$ are of the same degree $d$ having
rational roots $1-r_{\ell}$ and $1-q_{m}$ respectively with $0<r_{\ell}%
,q_{m}\leq1$ and $r_{\ell}\neq q_{m}$ for $\ell,m=1,2,...,d,$ so that
\{$r_{1}$,$r_{2}$...,$r_{d}$\}\ and\textsf{\ }\{$q_{1}$,$q_{2}$...,$q_{d}$\}
are disjoint rational multisets where, if a denominator $b$ occurs, all
irreducible fractions $a/b\in(0,1]$ also occur (i.e. either $a=b=1$ or $0<a<b$
and $\gcd(a,b)=1$). With this setup a given real constant $\omega
=\lim_{N\rightarrow\infty}\omega_{N},$ can be represented by the following
rational approximants sequence of normalized hypergeometric-type convergent
sums $\omega_{N}$\vspace*{-4pt}%
\begin{equation}
\omega_{N}=\frac{1}{\beta}\cdot%
{\displaystyle\sum\limits_{n=1}^{N}}
\mathcal{R(}n\mathcal{)\cdot H(}n)\text{ }\mathcal{=}\text{ }\frac{1}{\beta
}\cdot%
{\displaystyle\sum\limits_{n=1}^{N}}
\mathcal{R(}n\mathcal{)\cdot\rho}^{n}\,\mathcal{M}(n),\text{
\ \ {\small \textit{N}}}\in%
\mathbb{N}%
\label{1}%
\end{equation}
where $\beta\in%
\mathbb{Z}
_{\neq0}$ is a searching normalization factor, $\mathcal{R(}n\mathcal{)}%
=p(n)/r(n)$ is the rational part\ and $\mathcal{H(}n\mathcal{)}=\mathcal{\rho
}^{n}\,\mathcal{M}(n)$ is the hypergeometric motive factor where $|\rho
|\leq1,$ $\rho$ $\in%
\mathbb{Q}
$ is the convergence rate and $\mathcal{M(}n\mathcal{)}$ is defined by the
following ratio of products of $d$ Pochhammer's symbols$\vspace*{-2pt}$%
\begin{equation}
\mathcal{M}(n)=\left[
\begin{array}
[c]{cccc}%
r_{1} & r_{2} & ... & r_{d}\\
q_{1} & q_{2} & ... & q_{d}%
\end{array}
\vspace*{1pt}\right]  _{n}=\dfrac{(r_{1})_{n}\,(r_{2})_{n}\,...\,(r_{d})_{n}%
}{(q_{1})_{n}\,(q_{2})_{n}\,...\,(q_{d})_{n}}\vspace*{6pt}\label{2}%
\end{equation}
where $d$ is the hypergeometric deepness and $(\upsilon)_{n}=\upsilon
(\upsilon+1)\ldots(\upsilon+n-1)$ is the rising factorial. This means
that\ there exist leading coefficients $c_{r},c_{q}\in%
\mathbb{Z}
_{\neq0}$ whose ratio holds $c_{r}/c_{q}=\rho=\lim_{n\rightarrow\infty
}r(n)/q(n)$ and\vspace*{-4pt}%

\begin{equation}%
\begin{array}
[c]{ccc}%
r(n) & = & c_{r}(n-1+r_{1})(n-1+r_{2})\ldots(n-1+r_{d})%
\genfrac{}{}{0pt}{0}{\genfrac{}{}{0pt}{1}{{}}{{}}}{{}}%
\\
q(n) & = & c_{q}(n-1+q_{1})(n-1+q_{2})\ldots(n-1+q_{d})%
\genfrac{}{}{0pt}{0}{\genfrac{}{}{0pt}{1}{{}}{{}}}{{}}%
\end{array}
\vspace*{1pt}\label{2a}%
\end{equation}
Note that rational roots $1-r_{\ell}$ and $1-q_{m}$ enter as elements
$r_{\ell}$ and $q_{m}$ in $\mathcal{M}(n)$, therefore$\vspace*{-2pt}$%
\begin{equation}
\mathcal{H(}n\mathcal{)}=\rho^{n}\cdot\dfrac{(r_{1})_{n}\,(r_{2}%
)_{n}\,...\,(r_{d})_{n}}{(q_{1})_{n}\,(q_{2})_{n}\,...\,(q_{d})_{n}}=%
{\displaystyle\prod\limits_{k=1}^{n}}
\frac{r(k)}{q(k)}\vspace*{0pt}\label{2b}%
\end{equation}
so the hypergeometric-type series can be equivalently written either
as$\vspace*{-2pt}$%
\begin{equation}
\omega=\frac{1}{\beta}\cdot\lim_{N\rightarrow\infty}%
{\displaystyle\sum\limits_{n=1}^{N}}
\frac{p(n)}{r(n)}\cdot\mathcal{\rho}^{n}\mathcal{\cdot}\left[
\begin{array}
[c]{cccc}%
r_{1} & r_{2} & ... & r_{d}\\
q_{1} & q_{2} & ... & q_{d}%
\end{array}
\vspace*{1pt}\right]  _{n}\vspace*{-3pt}\label{7a}%
\end{equation}
or\vspace{-6pt}%
\begin{equation}
\omega=\frac{1}{\beta}\cdot\lim_{N\rightarrow\infty}%
{\displaystyle\sum\limits_{n=1}^{N}}
\frac{p(n)}{r(n)}\cdot%
{\displaystyle\prod\limits_{k=1}^{n}}
\frac{r(k)}{q(k)}\vspace*{3pt}\label{7b}%
\end{equation}
This last one is the classical form to apply the Binary Splitting algorithm
\cite{PAPA} and variations \cite{CHENG} where the computation time to evaluate
constant $\omega$ with $B$ bit is $\mathcal{O}(\mathit{M}(B)\cdot\log B)$,
being \textit{M}$(B)$ the time needed to multiply two $B$-bit numbers. A
practical measurement associated to \textit{M}$(B)$, that allows to
(asymptotically) rank, classify and compare different hypergeometric-type
algorithms by performance is\ given by the \textit{Binary Splitting Cost
}\cite{AYEE}%
\begin{equation}
C_{s}=-\frac{4\cdot d}{\log|\rho|}\label{3}%
\end{equation}
$C_{s}$ provides a method to know a priori which algorithm (formula, identity)
is better or faster than another, being a fundamental quantity to search, find
and select the most efficient hypergeometric series to calculate a given
constant $\omega$.

\section{Searching Process and Integer Relationships Detection\medskip}

The searching strategy starts from $\mathcal{M}(n)$ Eq.(\ref{2}), with a
selected sequence of small hypergeometric deepness $d=2,3,4,5,6$ that defines
the cardinality of the hypergeometric parameters multisets. For each size $d,$
polynomials $r(n)$ and $q(n)$ are built by specifying their hypergeometric
motive cyclotomic parameters [\textit{hgmcyclo} PARI function \cite{PARI}]
that expand into blocks containing all irreducible rational root elements
$1-t$ (blocks entering in $\mathcal{M}(n)$ as $t$). All cyclotomic parameters
for that $d$ are iterated by partitioning until all partitions are
exhausted$.$ For each parameters partition the rational convergence rate is
built as $\rho=\pm\,p_{1}^{n_{1}}\cdot p_{2}^{n_{2}}\cdot...\cdot\,p_{\ell
}^{n_{\ell}}$ where $\{p_{i}\}_{i=1}^{\ell}$ are primes extracted from all the
denominators of the irreducible rational roots and the powers $\{n_{1}%
,n_{2},...,n_{\ell}\}\in%
\mathbb{Z}
^{\ell}$ form the lattice iterating $\rho$ by brute force fulfilling
$n_{i\text{\thinspace min}}\leq$ $n_{i}\leq n_{i\text{\thinspace max}}$,
$C_{s}<$ $C_{0}$ Eq.(\ref{3}), $N<N_{0}$ Eq.(\ref{3b}) and $|\rho|<\rho_{0}$
for given exponent limits $n_{i\text{\thinspace min}},$ $n_{i\text{\thinspace
max}}$\ and stated bounds $C_{0},$ $N_{0},$ $\rho_{0}$, so\vspace*{-6pt}%
\begin{equation}
\log|\rho|=%
{\displaystyle\sum\limits_{i=0}^{\ell}}
n_{i}\cdot\log p_{i}<L\vspace*{-4pt}\label{3z}%
\end{equation}
where $L=\min(\log\rho_{0},-4d/C_{0},-(h+2)\cdot w_{p}\cdot\log$%
{\small (10)}$/N_{0})$ and $(h+2)\cdot w_{p}$ are the decimal digits of
working precision needed (see below). Note that for logarithms $\omega=\log
p,$ prime factors of $p-1,$ $p$ and $p+1$ must be also included in
$\{p_{i}\}_{i=1}^{\ell}$--see Eq.(\ref{7}) and Eqs.(\ref{18z}--\ref{18y}%
)--.\thinspace\ Finally, the search terminates if, for fixed $d,q(n),r(n),\rho
,$ an integer $\beta$ and a polynomial $p(n)$ of degree $h$
\begin{equation}
p(n)=\alpha_{0}+\alpha_{1}n+\alpha_{2}n^{2}+...+\alpha_{h}n^{h},\text{
\ \ }\alpha_{h}\neq0\vspace*{4pt}\label{3a}%
\end{equation}
having integer coefficients $\alpha_{0},\alpha_{1},...\alpha_{h},$ are found,
so that Eq.(\ref{7b}) can be expressed as\vspace*{-4pt}%
\begin{equation}%
{\displaystyle\sum\limits_{i=0}^{h}}
\alpha_{i}\,s_{i}-\beta\,\omega=0\label{3x}%
\end{equation}
that represents an integer relationship, if it exists, between the real
constant $\omega$ and the real values $s_{i},$ $i=0,1,2,...,h$ given by the
sums\vspace*{-4pt}
\begin{equation}
s_{i}=\lim_{N\rightarrow\infty}%
{\displaystyle\sum\limits_{n=1}^{N}}
\frac{n^{i}}{r(n)}\cdot%
{\displaystyle\prod\limits_{k=1}^{n}}
\frac{r(k)}{q(k)}\label{3b}%
\end{equation}
There are $h+2$ integer values $\beta,\alpha_{0},\alpha_{1},...,\alpha_{h}$ to
be detected. For a given deepness $d$ and a constant $\omega$ satisfying
Eq.(\ref{7a}) the value of $h$ is obtained as $h=d-w_{t}$ where\ $w_{t}$ is
the known constant's weight or period's degree, (H. Cohen\thinspace
\cite{HCOHEN}, Def. 1.4). So weight $w_{t}=0$ are all algebraic numbers;
$w_{t}=1$, $\pi$, lemniscate $\varpi$, logarithms of algebraic numbers;
$w_{t}=2$, Catalan's constant $G$, $\Gamma(1/3)^{3}$, Dirichlet's constants
$L(-3,2)$, $L(-7,2)$; $w_{t}=3$, $\zeta(3)$, $\Gamma(1/4)^{4}$; $w_{t}=$
$n,$\ $\zeta(n)$, polylogarithms $Li_{n}(1/m)\ $with $m $ integer and so on.
For logarithms $\omega=\log p$ with $p\in%
\mathbb{Z}
_{>1}$ the polynomial degree\ is $h=d-1$.

\bigskip

There are some integer relations detection algorithms, PSLQ \cite{PSLQ}%
\cite{FERGUSON}, HJLS \cite{HJLS} and LLL (lattice reduction) \cite{LLL},
among others, that can be used to search solutions to Eq.(\ref{3x}). All of
them need at least a rough working precision $w_{p}=\log_{10}(\max
(|\beta|,|\alpha_{0}|,|\alpha_{1}|,...,|\alpha_{h}|))$ per integer to be
detected. So for $64$ bit integer coefficients, a decimal digits working
precision higher than $(h+2)\cdot w_{p}=\,(h+2)\cdot\log_{10}2^{64}%
\,\simeq\,20\cdot(d-w_{t}+2)$ is required. A custom-made very efficient script (\texttt{myPSLQ}, available as it is indicated in Section 8.)
for searching hypergeometric-type identities was coded in PARI--GP
\cite{PARI}, a specialized Number Theory platform that works fine with high
precision computing and has a very fast native implementation of LLL
[\textit{lindep} function], the lattice reduction algorithm.

\section{A variation of the Wilf--Zeilberger method\medskip}

A different approach can be taken to obtain new proven identities for constant
$\omega$ through the application of an extension of the Wilf--Zeilberger
method beyond \cite{ZEIL}. \ For a generic $\omega$ and a known base WZ pair
$(F(n,k),G(n,k))$ a family of WZ pairs $(F_{s,t}(n,k),G_{s,t}(n,k))$ with
$s\in%
\mathbb{N}
$ and $t\in%
\mathbb{Z}
$ is built fulfilling these limit identities whenever the flawless series
\cite{FLAW} converge\vspace*{-4pt}%
\begin{equation}%
\begin{tabular}
[c]{lllll}%
$\lim\limits_{n\rightarrow\infty}%
{\displaystyle\sum\limits_{k=0}^{n}}
F(n,k)$ & $=$ & $\lim\limits_{n\rightarrow\infty}%
{\displaystyle\sum\limits_{k=0}^{n}}
F_{s,t}(n,k)$ & $=$ & $0\vspace*{4pt}$\\
$\lim\limits_{k\rightarrow\infty}%
{\displaystyle\sum\limits_{n=0}^{k}}
G(n,k)$ & $=$ & $\lim\limits_{k\rightarrow\infty}%
{\displaystyle\sum\limits_{n=0}^{k}}
G_{s,t}(n,k)$ & $=$ & $0%
\genfrac{}{}{0pt}{0}{\genfrac{}{}{0pt}{1}{{}}{{}}}{{}}%
\vspace*{-4pt}$%
\end{tabular}
\label{4a}%
\end{equation}

\begin{equation}
\omega\ \,=\,\
{\displaystyle\sum\limits_{k=0}^{\infty}}
F(0,k)\ \,=\,\
{\displaystyle\sum\limits_{k=0}^{\infty}}
F_{s,t}(0,k)\,\ =\,\
{\displaystyle\sum\limits_{n=0}^{\infty}}
G_{s,t}(n,0)\,\ =\,\
{\displaystyle\sum\limits_{n=0}^{\infty}}
G(n,0)\vspace*{-2pt}\label{4}%
\end{equation}
By selecting
\begin{equation}
F_{s,t}(n,k)\ =\ F(s\hspace*{1pt}n,k+t\hspace*{1pt}n)\vspace*{2pt}\label{5}%
\end{equation}
these identities are found\vspace*{-2pt}
\begin{equation}%
\begin{tabular}
[c]{lllll}%
$G_{s,t}(n,k)$ & $=$ & $%
{\displaystyle\sum\limits_{\ell=0}^{s-1}}
G(s\hspace*{1pt}n+\ell,k+t\hspace*{1pt}n)$ & $+$ & $%
{\displaystyle\sum\limits_{m=0}^{t-1}}
F(s\hspace*{1pt}n+s,m+k+t\hspace*{1pt}n)\vspace*{3pt}$\\
& $=$ & $%
{\displaystyle\sum\limits_{\ell=0}^{s-1}}
G(s\hspace*{1pt}n+\ell,k+t+t\hspace*{1pt}n)$ & $+$ & $%
{\displaystyle\sum\limits_{m=0}^{t-1}}
F(s\hspace*{1pt}n,m+k+t\hspace*{1pt}n)$%
\end{tabular}
\label{5a}%
\end{equation}
This $(s,t)${\small --}$W\hspace*{-2pt}Z$ pair transformation can be written
as $W\hspace*{-2pt}Z_{s,t}:$ $(F,G)\rightarrow(F_{s,t},G_{s,t})$. Note
that\ the discovered relationships $\omega=%
{\textstyle\sum\nolimits_{n=0}^{\infty}}
G_{s,t}(n,0)$ Eq.(\ref{4}), always result proven since $W\hspace*{-2pt}%
Z_{s,t}$ provides the corresponding rational certificates straightwise. There
are some techniques to guess a proper bivariate companion $F(n,k)$ from a
known $\omega=%
{\textstyle\sum\nolimits_{k=0}^{\infty}}
F(0,k)$ identity such that a WZ pair $(F,G)$ exists (for instance, Au's WZ
seeds \cite{AU}) which allows to apply the $W\hspace*{-2pt}Z_{s,t}$
transformation to get new families of identities. Everything written in the
previous sections and in this one applies to generic constants $\omega$.
Specializing it for logarithms $\omega=\log p,$ a companion $F(n,k)$ can be
raised from these two linearly\thinspace convergent Gauss hypergeometric
series linked by\thinspace${}_{2}F_{1}$\thinspace Pfaff's
transformations\thinspace\cite{EULER}\vspace*{2pt}
\begin{equation}%
\begin{tabular}
[c]{rrlll}%
$\log p$ & $=$ & {\small 2 }$\cdot\frac{p-1}{p+1}\cdot{}_{2}F_{1}\left(
\frac{1}{2},1;\frac{3}{2};(\frac{p-1}{p+1})^{2}\right)  \vspace*{6pt}$ & $%
\genfrac{}{}{0pt}{0}{\genfrac{}{}{0pt}{1}{{}}{{}}}{{}}%
$ & $p>0$\\
$\log p$ & $=$ & $\frac{p^{2}-1}{2p}\cdot{}_{2}F_{1}\left(  1,1;\frac{3}%
{2};-\frac{(p-1)^{2}}{4p}\right)  $ & $%
\genfrac{}{}{0pt}{0}{\genfrac{}{}{0pt}{1}{{}}{{}}}{{}}%
\vspace*{2pt}$ & $|p-3|<2\sqrt{2}$%
\end{tabular}
\label{6}%
\end{equation}
whose Taylor's terms can be associated to $F(0,k)$ and $G(n,0)$ respectively
(or viceversa)$\vspace*{-1pt}.$ Both series are easily\vspace*{-2pt}
interlaced through the beta function $B(u,v)=\frac{\Gamma(u)\Gamma(v)}%
{\Gamma(u+v)}$ to lift up two\ different $F(n,k)$ giving $(F,G)$ pairs that
hold Eq.(\ref{4a}) and work for\ $|p-3|<2\sqrt{2}$ as
\begin{equation}%
\begin{tabular}
[c]{rrl}%
$F^{\text{{\tiny (1)}}}(n,k)$ & $=$ & {\small \ }$(-1)^{n}\cdot\dfrac
{(p-1)^{2n+2k+1}}{(4p)^{n}\,(p+1)^{2k+1}}\cdot B(k+\tfrac{1}{2},n+1)\vspace
*{8pt}$\\
$F^{\text{{\tiny (2)}}}(n,k)$ & $=$ & $(-1)^{k}\cdot\dfrac{(p-1)^{2n+2k+1}%
}{(4p)^{k+1}\,(p+1)^{2n-1}}\cdot{}B(k+1,n+\tfrac{1}{2})$%
\end{tabular}
\label{7}%
\end{equation}
Due to convergence of the source series only $p=2,3,4$ and $5$ are valid for
$p\in%
\mathbb{Z}
_{>1}.$\vspace*{2pt} The resulting\ $W\hspace*{-2pt}Z_{s,t}$
transformations\ will be referred either as $\left.  W\hspace*{-2pt}%
Z_{s,t}^{\text{{\tiny (1)}}}\right\vert _{p}$ or $\left.  W\hspace
*{-2pt}Z_{s,t}^{\text{{\tiny (2)}}}\right\vert _{p}$ respectively.

\section{log$\hspace{2pt}p$ fast series, $p\ =$ 2, 3, 5\medskip}

From the huge number of hypergeometric-type logarithm series found only those
that perform faster than Machin-type algorithms and have the lowest binary
splitting costs are selected. Three series with $d=2$ for $\log2,$ $\log3$ and
$\log5$, four identities with $d=4$ for $\log2$ and $\log3$ and one with $d=6$
for $\log2.$ Except for one, all of them are proven. The proofs of the
identities with $d=2$ are placed in \textit{Section} 6, further on.\thinspace
They are based on the beta integral (BI) \cite{GUIGOU}, hypergeometric closed
forms (CF) \cite{SUNZHOU}\cite{CHBROWN} and on Wilf--Zeilberger $W\hspace
*{-2pt}Z_{s,t}$ certificates as well. The proofs for $d=4,$ $6$ are obtained
by BI and also by $W\hspace*{-2pt}Z_{s,t}$ certificates except one with $d=4$
that remains unproven.

\subsection{$\boldsymbol{d=2}$}

These are three Ramanujan-type identities that were found for first time on
Dec. 2023 by applying the searching algorithm, all of them belong to the same
family,
\begin{equation}%
\begin{tabular}
[c]{rrr}%
$\log2$ & $=$ & $%
{\displaystyle\sum\limits_{n=1}^{\infty}}
\dfrac{1794\,n-297}{2n\,(2n-1)}\cdot\left(  \dfrac{1}{2^{4}\cdot3^{5}}\right)
^{n}\cdot%
\begin{bmatrix}%
\genfrac{}{}{0pt}{0}{{}}{{}}%
1 & \frac{1}{2}%
\genfrac{}{}{0pt}{0}{{}}{{}}%
\\%
\genfrac{}{}{0pt}{0}{{}}{{}}%
\frac{1}{6} & \frac{5}{6}%
\genfrac{}{}{0pt}{0}{{}}{{}}%
\end{bmatrix}
_{n}$%
\end{tabular}
\label{8}%
\end{equation}
This series has a cost $C_{s}=0.9679\ldots$ being the fastest to compute this
constant. It has been incorporated into y-cruncher software \cite{YCRUNCHER}
as the primary algorithm for $\log2$ and it was employed to get $3\cdot
10^{12}$ decimal places on Feb. 2024 \cite{NWORLD}. (Part of this story is
told in \cite{ZUN2}). This identity has also been included in FLINT
\cite{FLINT} (Release 3.3.0-dev. See 9.31.2, p.774. Jan 2025).\vspace
*{0pt}\ Eq.(\ref{8}) is proven by BI, \vspace*{0pt}CF methods and also
by\ $W\hspace*{-2pt}Z_{2,1}^{\text{{\tiny (1)}}}|_{p=2}$, $W\hspace
*{-2pt}Z_{1,2}^{\text{{\tiny (2)}}}|_{p=2}$ certificates.\vspace*{-2pt}%
\begin{equation}%
\begin{tabular}
[c]{rrr}%
$\log3$ & $=$ & $%
{\displaystyle\sum\limits_{n=1}^{\infty}}
\dfrac{88\,n-14}{n\,(2n-1)}\cdot\left(  \dfrac{1}{3^{5}}\right)  ^{n}\cdot%
\begin{bmatrix}%
\genfrac{}{}{0pt}{0}{{}}{{}}%
1 & \frac{1}{2}%
\genfrac{}{}{0pt}{0}{{}}{{}}%
\\%
\genfrac{}{}{0pt}{0}{{}}{{}}%
\frac{1}{6} & \frac{5}{6}%
\genfrac{}{}{0pt}{0}{{}}{{}}%
\end{bmatrix}
_{n}$%
\end{tabular}
\label{8a}%
\end{equation}
This identity has a cost $C_{s}=1.4564\ldots$ being currently the fastest
single series for $\log3$. It has been included into y-cruncher software
\cite{YCRUNCHER} \vspace*{0pt}as the main (primary) algorithm.\ This identity
is proven by BI, CF methods and by\ $W\hspace*{-2pt}Z_{2,1}^{\text{{\tiny (1)}%
}}|_{p=3}$, $W\hspace*{-2pt}Z_{1,2}^{\text{{\tiny (2)}}}|_{p=3}$ certificates
as well.%
\begin{equation}%
\begin{tabular}
[c]{rrr}%
$\log5$ & $=$ & $%
{\displaystyle\sum\limits_{n=1}^{\infty}}
\dfrac{364\,n-62}{-n\,(2n-1)}\cdot\left(  \dfrac{-1}{3^{3}\cdot5^{2}}\right)
^{n}\cdot%
\begin{bmatrix}%
\genfrac{}{}{0pt}{0}{{}}{{}}%
1 & \frac{1}{2}%
\genfrac{}{}{0pt}{0}{{}}{{}}%
\\%
\genfrac{}{}{0pt}{0}{{}}{{}}%
\frac{1}{6} & \frac{5}{6}%
\genfrac{}{}{0pt}{0}{{}}{{}}%
\end{bmatrix}
_{n}$%
\end{tabular}
\label{8b}%
\end{equation}
This identity has a cost $C_{s}=1.2280\ldots$ and it is the most efficient
single series for $\log5$ (Jan. 2025). It has been included into y-cruncher
software \cite{YCRUNCHER} as the primary algorithm to compute this
constant.\ \vspace*{-1pt}This series is proven by BI, CF methods and also by
$W\hspace*{-2pt}Z_{2,1}^{\text{{\tiny (1)}}}|_{p=2\pm i}$, $W\hspace
*{-2pt}Z_{1,2}^{\text{{\tiny (2)}}}|_{p=2\pm i}$ certificates.\ (See Sections
6.3, Alternating Series and 6.4 further on)

\subsection{$\boldsymbol{d=4}$}

The following identities were mainly found (and proven) in Nov. 2024 by means
of $W\hspace*{-2pt}Z_{s,t}$ transformations from Eq.(\ref{7}) and the
application of a special Maplesoft$\texttrademark$ Maple \cite{MAPLE} code
using Briochet--Salvy implementation of Creative Telescoping \cite{CRTEL}%
\begin{equation}%
\begin{tabular}
[c]{r}%
$\log2\,=\,\,%
{\displaystyle\sum\limits_{n=1}^{\infty}}
\dfrac{P(n)}{2n(2n-1)(6n-1)(6n-5)}\cdot\left(  \dfrac{1}{2^{4}\cdot3^{3}%
\cdot5^{5}}\right)  ^{n}\,%
\begin{bmatrix}%
\genfrac{}{}{0pt}{0}{{}}{{}}%
1 & \frac{1}{2} & \frac{1}{6} & \frac{5}{6}%
\genfrac{}{}{0pt}{0}{{}}{{}}%
\\%
\genfrac{}{}{0pt}{0}{{}}{{}}%
\frac{1}{10} & \frac{3}{10} & \frac{7}{10} & \frac{9}{10}%
\genfrac{}{}{0pt}{0}{{}}{{}}%
\end{bmatrix}
_{n}$%
\end{tabular}
\label{9}%
\end{equation}%
\[%
\begin{tabular}
[c]{lll}%
$P(n)$ & $=$ & $13885704\,n^{3}-15397068\,n^{2}+4353342\,n-295245%
\genfrac{}{}{0pt}{0}{\genfrac{}{}{0pt}{1}{{}}{{}}}{{}}%
$%
\end{tabular}
\]
Eq.(\ref{9}) has $C_{s}=1.1335\ldots$ being currently --Jan. 2025-- the second
fastest single series, just after Eq(\ref{8}), to compute this constant.

Eq.(\ref{9}) is included in y-cruncher v0.8.6 \cite{YCRUNCHER} as the
secondary algorithm for digits validation of $\log2.$ \vspace*{0pt}This series
is found and proven applying $W\hspace*{-2pt}Z_{2,3}^{\text{{\tiny (1)}}%
}|_{p=2}$ or $W\hspace*{-2pt}Z_{3,2}^{\text{{\tiny (2)}}}|_{p=2}$%
\vspace*{-6pt}%

\begin{equation}%
\begin{tabular}
[c]{r}%
$\log2\,=\,%
{\displaystyle\sum\limits_{n=1}^{\infty}}
\dfrac{P(n)}{4n(2n-1)(4n-1)(4n-3)}\cdot\left(  \dfrac{1}{2^{4}\cdot3^{2}%
\cdot5^{5}}\right)  ^{n}\,%
\begin{bmatrix}%
\genfrac{}{}{0pt}{0}{{}}{{}}%
1 & \frac{1}{2} & \frac{1}{4} & \frac{3}{4}%
\genfrac{}{}{0pt}{0}{{}}{{}}%
\\%
\genfrac{}{}{0pt}{0}{{}}{{}}%
\frac{1}{10} & \frac{3}{10} & \frac{7}{10} & \frac{9}{10}%
\genfrac{}{}{0pt}{0}{{}}{{}}%
\end{bmatrix}
_{n}$%
\end{tabular}
\label{11}%
\end{equation}%
\[%
\begin{tabular}
[c]{lll}%
$P(n)$ & $=$ & $3927264\,n^{3}-4300512\,n^{2}+1209726\,n-81891%
\genfrac{}{}{0pt}{0}{\genfrac{}{}{0pt}{1}{{}}{{}}}{{}}%
$%
\end{tabular}
\]
Eq.(\ref{11}) is a fast series for $\log2$ having a very low cost
$C_{s}=1.2292\ldots$ This identity is found and proven through $W\hspace
*{-2pt}Z_{4,1}^{\text{{\tiny (1)}}}|_{p=2}$ or $W\hspace*{-2pt}Z_{1,4}%
^{\text{{\tiny (2)}}}|_{p=2}$%

\begin{equation}%
\begin{tabular}
[c]{r}%
$\log2\,=\,%
{\displaystyle\sum\limits_{n=1}^{\infty}}
\dfrac{P(n)}{3n(2n-1)(3n-1)(3n-2)}\cdot\left(  \dfrac{1}{2^{13}\cdot3^{3}%
}\right)  ^{n}\cdot%
\begin{bmatrix}%
\genfrac{}{}{0pt}{0}{{}}{{}}%
1 & \frac{1}{2} & \frac{1}{3} & \frac{2}{3}%
\genfrac{}{}{0pt}{0}{{}}{{}}%
\\%
\genfrac{}{}{0pt}{0}{{}}{{}}%
\frac{1}{12} & \frac{5}{12} & \frac{7}{12} & \frac{11}{12}%
\genfrac{}{}{0pt}{0}{{}}{{}}%
\end{bmatrix}
_{n}$%
\end{tabular}
\label{13}%
\end{equation}%
\[%
\begin{tabular}
[c]{lll}%
$P(n)$ & $=$ & $686430\,n^{3}-742257\,n^{2}+223397\,n-13858%
\genfrac{}{}{0pt}{0}{\genfrac{}{}{0pt}{1}{{}}{{}}}{{}}%
$%
\end{tabular}
\]
Identity Eq.(\ref{13}) is a fast series found on Dec. 2023 applying the
searching process via LLL algorithm. It has a very low cost $C_{s}%
=1.3001\ldots$ This was the secondary digits validation algorithm for $\log2$
that was included in native form inside y-cruncher v0.8.5 (2024)
\cite{YCRUNCHER}, but now it has been superseeded by Eq.(\ref{9}).
Eq.(\ref{13}) was applied on Feb. 2024 to verify 3 trillions of decimal places
record for this constant. This identity is just conjectured since it has not
been possible to prove it so far despite some attempts \cite{ZUN1}.%

\begin{equation}%
\begin{tabular}
[c]{r}%
$\log3\,=\,%
{\displaystyle\sum\limits_{n=1}^{\infty}}
\dfrac{P(n)}{n(2n-1)(6n-1)(6n-5)}\cdot\left(  \dfrac{3}{2^{4}\cdot5^{5}%
}\right)  ^{n}\cdot%
\begin{bmatrix}%
\genfrac{}{}{0pt}{0}{{}}{{}}%
1 & \frac{1}{2} & \frac{1}{6} & \frac{5}{6}%
\genfrac{}{}{0pt}{0}{{}}{{}}%
\\%
\genfrac{}{}{0pt}{0}{{}}{{}}%
\frac{1}{10} & \frac{3}{10} & \frac{7}{10} & \frac{9}{10}%
\genfrac{}{}{0pt}{0}{{}}{{}}%
\end{bmatrix}
_{n}$%
\end{tabular}
\label{15a}%
\end{equation}%
\[%
\begin{tabular}
[c]{lll}%
$P(n)$ & $=$ & $141168\,n^{3}-158016\,n^{2}+44804\,n-3040%
\genfrac{}{}{0pt}{0}{\genfrac{}{}{0pt}{1}{{}}{{}}}{{}}%
$%
\end{tabular}
\]
Series Eq.(\ref{15a}) has $C_{s}=1.6459\ldots$ being currently the second
fastest, just after Eq.(\ref{8a}), to compute $\log3$. It has been included
natively in y-cruncher v0.8.6 (2025) \cite{YCRUNCHER} as the secondary
algorithm for digits validation of this constant. It performs faster than any
Machin-type formula to compute $\log3.$\vspace*{-1pt} It is found and proven
by $W\hspace*{-2pt}Z_{2,3}^{\text{{\tiny (1)}}}|_{p=3}$ or $W\hspace
*{-2pt}Z_{3,2}^{\text{{\tiny (2)}}}|_{p=3}.$

\subsection{$\boldsymbol{d=6}$}

The following series has a very low cost $C_{s}=1.2189\ldots$ It was found and
proven from Eq.(\ref{7}) by means of $W\hspace*{-2pt}Z_{4,3}%
^{\text{{\tiny (1)}}}|_{p=2}$ or $W\hspace*{-2pt}Z_{3,4}^{\text{{\tiny (2)}}%
}|_{p=2}.\vspace{-1pt}$ (See also Subsection 6.5 Example I, for a BI proof). 
This is the 3$^{rd}$ fastest single series formula to compute
$\log2$ so far$.$
\begin{equation}%
\begin{tabular}
[c]{rrr}%
$\log2$ & $=$ & $\dfrac{1}{4}\cdot%
{\displaystyle\sum\limits_{n=1}^{\infty}}
\dfrac{P(n)}{R(n)}\cdot\left(  \dfrac{1}{2^{4}\cdot3^{3}\cdot7^{7}}\right)
^{n}\cdot%
\begin{bmatrix}%
\genfrac{}{}{0pt}{0}{{}}{{}}%
1 & \frac{1}{2} & \frac{1}{4} & \frac{3}{4} & \frac{1}{6} & \frac{5}{6}%
\genfrac{}{}{0pt}{0}{{}}{{}}%
\\%
\genfrac{}{}{0pt}{0}{{}}{{}}%
\frac{1}{14} & \frac{3}{14} & \frac{5}{14} & \frac{9}{14} & \frac{11}{14} &
\frac{13}{14}%
\genfrac{}{}{0pt}{0}{{}}{{}}%
\end{bmatrix}
_{n}$%
\end{tabular}
\label{18}%
\end{equation}
where polynomials $P(n)$ and $R(n)$ are%
\[%
\begin{tabular}
[c]{lll}%
$P(n)$ & $=$ & $81969540480\,n^{5}-169950180480\,n^{4}+126495134424\,n^{3}%
\genfrac{}{}{0pt}{0}{\genfrac{}{}{0pt}{1}{{}}{{}}}{{}}%
$\\
&  & $-40884797604\,n^{2}+5510613042\,n-226846575%
\genfrac{}{}{0pt}{0}{\genfrac{}{}{0pt}{1}{{}}{{}}}{{}}
$\\
$R(n)$ & $=$ & $n(2n-1)(4n-1)(4n-3)(6n-1)(6n-5)%
\genfrac{}{}{0pt}{0}{\genfrac{}{}{0pt}{1}{{}}{{}}}{{}}%
$%
\end{tabular}
\]

\section{Proofs for $d=2$\bigskip}

\begin{center}%
\begin{tabular}
[c]{cccccccccccc}\hline\hline
\multicolumn{1}{||c}{\textit{p}} & \multicolumn{1}{|c}{$\alpha$} &
\multicolumn{1}{|c}{$\beta$} & \multicolumn{1}{|c}{$\gamma$} &
\multicolumn{1}{|c}{\textit{a}} & \multicolumn{1}{|c}{\textit{b}} &
\multicolumn{1}{|c}{\textit{c}} & \multicolumn{1}{|c}{{\small \textit{A}}%
$_{1}$} & \multicolumn{1}{|c}{{\small \textit{A}}$_{2}$} &
\multicolumn{1}{|c}{{\small \textit{A}}$_{3}$} & \multicolumn{1}{|c}{$\rho$} &
\multicolumn{1}{|c||}{\textit{z}$%
\genfrac{}{}{0pt}{0}{\genfrac{}{}{0pt}{1}{{}}{{}}}{{}}%
$}\\\hline\hline
\multicolumn{1}{||c}{{\small 2}} & \multicolumn{1}{|r|}{{\small 1794}} &
\multicolumn{1}{|r|}{{\small --297}} & \multicolumn{1}{|r}{{\small 2}} &
\multicolumn{1}{|r}{{\small 598}} & \multicolumn{1}{|r}{{\small 499}} &
\multicolumn{1}{|r}{{\small 144}} & \multicolumn{1}{|c}{$\frac{25}{72}$} &
\multicolumn{1}{|r}{$-\frac{1}{192}$} & \multicolumn{1}{|r}{$\frac{1}{192}$} &
\multicolumn{1}{|c}{$\frac{1}{3888}$} & \multicolumn{1}{|c||}{{\small 3}$%
\genfrac{}{}{0pt}{0}{\genfrac{}{}{0pt}{1}{{}}{{}}}{{}}%
$}\\
\multicolumn{1}{||c}{{\small 3}} & \multicolumn{1}{|r|}{{\small 88}} &
\multicolumn{1}{|r|}{{\small --14}} & \multicolumn{1}{|r}{{\small 1}} &
\multicolumn{1}{|r}{{\small 176}} & \multicolumn{1}{|r}{{\small 148}} &
\multicolumn{1}{|r}{{\small 27}} & \multicolumn{1}{|c}{$\frac{5}{9}$} &
\multicolumn{1}{|r}{$-\frac{1}{18}$} & \multicolumn{1}{|r}{$\frac{1}{18}$} &
\multicolumn{1}{|c}{$\frac{1}{243}$} & \multicolumn{1}{|c||}{{\small 2}$%
\genfrac{}{}{0pt}{0}{\genfrac{}{}{0pt}{1}{{}}{{}}}{{}}%
$}\\
\multicolumn{1}{||c}{{\small 5}} & \multicolumn{1}{|r}{{\small --364}} &
\multicolumn{1}{|r}{{\small 62}} & \multicolumn{1}{|r}{{\small 1}} &
\multicolumn{1}{|r}{{\small 728}} & \multicolumn{1}{|r}{{\small 604}} &
\multicolumn{1}{|r}{{\small 75}} & \multicolumn{1}{|c}{$\frac{4}{5}$} &
\multicolumn{1}{|r}{$\frac{1}{25}$} & \multicolumn{1}{|r}{$-\frac{1}{25}$} &
\multicolumn{1}{|c}{{\small --}$\frac{1}{675}$} &
\multicolumn{1}{|c||}{{\small 2}\textit{i}$%
\genfrac{}{}{0pt}{0}{\genfrac{}{}{0pt}{1}{{}}{{}}}{{}}%
$}\\
\multicolumn{1}{||c}{{\small 7}} & \multicolumn{1}{|r}{{\small 312}} &
\multicolumn{1}{|r}{{\small --16}} & \multicolumn{1}{|r}{{\small 81}} &
\multicolumn{1}{|r}{{\small 468}} & \multicolumn{1}{|r}{{\small 444}} &
\multicolumn{1}{|r}{{\small 49}} & \multicolumn{1}{|c}{$\frac{15}{14}$} &
\multicolumn{1}{|r}{$-\frac{243}{196}$} & \multicolumn{1}{|r}{$\frac{243}%
{196}$} & \multicolumn{1}{|c}{$\frac{27}{196}$} &
\multicolumn{1}{|c||}{{\small 4/3}$%
\genfrac{}{}{0pt}{0}{\genfrac{}{}{0pt}{1}{{}}{{}}}{{}}%
$}\\
\multicolumn{1}{||c}{{\small 10}} & \multicolumn{1}{|r|}{{\small --126}} &
\multicolumn{1}{|r|}{{\small 23}} & \multicolumn{1}{|r}{{\small 2}} &
\multicolumn{1}{|r}{{\small 1134}} & \multicolumn{1}{|r}{{\small 927}} &
\multicolumn{1}{|r}{{\small 80}} & \multicolumn{1}{|c}{$\frac{9}{8}$} &
\multicolumn{1}{|r}{$\frac{81}{320}$} & \multicolumn{1}{|r}{$-\frac{81}{320}$}
& \multicolumn{1}{|c}{{\small --}$\frac{1}{80}$} &
\multicolumn{1}{|c||}{\textit{i}$\sqrt{\text{{\small 5/3}}}%
\genfrac{}{}{0pt}{0}{\genfrac{}{}{0pt}{1}{{}}{{}}}{{}}%
$}\\\hline\hline
&  &  &  & \vspace*{-4pt} &  &  &  &  &  &  & \\
\multicolumn{12}{c}{\textit{Table }I\textit{. }$\log p,$\textit{\ }%
$d=2$\textit{\ Hypergeometric Series Parameters}}\\
\multicolumn{12}{c}{}%
\end{tabular}
\vspace*{-10pt}
\end{center}

The proofs of $d=2$ identities, Eqs.(\ref{8}--\ref{8a}--\ref{8b}), follow by
expressing them in these two equivalent forms\vspace*{-8pt}%
\begin{equation}%
\begin{tabular}
[c]{lll}%
$\log p$ & $=$ & $\dfrac{1}{\gamma}\cdot%
{\displaystyle\sum\limits_{n=1}^{\infty}}
\dfrac{\alpha n+\beta}{n(2n-1)}\cdot\rho^{n}\cdot%
\begin{bmatrix}%
\genfrac{}{}{0pt}{0}{{}}{{}}%
1 & \frac{1}{2}%
\genfrac{}{}{0pt}{0}{{}}{{}}%
\\%
\genfrac{}{}{0pt}{0}{{}}{{}}%
\frac{1}{6} & \frac{5}{6}%
\genfrac{}{}{0pt}{0}{{}}{{}}%
\end{bmatrix}
_{n}$%
\end{tabular}
\label{20}%
\end{equation}%
\begin{equation}%
\begin{tabular}
[c]{lll}%
$\log p$ & $=$ & $\dfrac{1}{c}\cdot%
{\displaystyle\sum\limits_{n=0}^{\infty}}
\dfrac{an+b}{(6n+1)(6n+5)}\cdot\rho^{n}\cdot%
\begin{bmatrix}%
\genfrac{}{}{0pt}{0}{{}}{{}}%
1 & \frac{1}{2}%
\genfrac{}{}{0pt}{0}{{}}{{}}%
\\%
\genfrac{}{}{0pt}{0}{{}}{{}}%
\frac{1}{6} & \frac{5}{6}%
\genfrac{}{}{0pt}{0}{{}}{{}}%
\end{bmatrix}
_{n}$%
\end{tabular}
\label{21}%
\end{equation}
where $\{\alpha,\beta,\gamma\}=v/\gcd(v)$ and $v=$ sgn$(\rho)$ $\cdot$
denom$(\rho)\cdot\{a,$ $b-a,$ $18\cdot c\cdot\rho\}$ with inverse
$\{a,b,c\}=u/\gcd(u)$ and $u=$ denom$(\rho)\cdot\{18\cdot\rho\cdot\alpha,$
$18\cdot\rho\cdot(\alpha+\beta),$ $\gamma\}.$ \textit{Table} I above contains
values of the parameters for these identities to be proven. Two slow
Ramanujan-type series for $\log7$ and $\log10$ that belong to the same family
and were found by applying the searching algorithm of Section 3, have also
been included.

\subsection{Beta Integral method}

Inspired on {\small Gour\'{e}vitch--Guillera }\cite{GUIGOU} beta
integral\ technique, name $\mathcal{S}$ the \textrm{RHS} of Eq.(\ref{21}) and
set the Pochhammer's ratio as\vspace*{4pt}
\begin{equation}%
\begin{bmatrix}%
\genfrac{}{}{0pt}{0}{{}}{{}}%
1 & \frac{1}{2}%
\genfrac{}{}{0pt}{0}{{}}{{}}%
\\%
\genfrac{}{}{0pt}{0}{{}}{{}}%
\frac{1}{6} & \frac{5}{6}%
\genfrac{}{}{0pt}{0}{{}}{{}}%
\end{bmatrix}
_{n}=\left(  \frac{27}{4}\right)  ^{n}\cdot\frac{\Gamma(2n+1)\hspace
*{1pt}\Gamma(n+\frac{1}{2})}{\Gamma(3n+\frac{1}{2})}\label{22}%
\end{equation}
the rational part
\begin{equation}
\mathcal{R(}n\mathcal{)={}}\frac{1}{c}\cdot\frac{an+b}{(6n+1)(6n+5)}\label{23}%
\end{equation}
can be written as%
\begin{equation}%
\begin{array}
[c]{c}%
\mathcal{R(}n\mathcal{)}=\dfrac{A_{1}}{3n+\frac{1}{2}}+\dfrac{A_{2}%
\cdot(2n+1)}{(3n+\frac{1}{2})(3n+\frac{3}{2})}+\dfrac{A_{3}\cdot
(2n+1)(2n+2)}{(3n+\frac{1}{2})(3n+\frac{3}{2})(3n+\frac{5}{2})}%
\end{array}
\label{23b}%
\end{equation}
where $A_{1},A_{2}$ and $A_{3}$ are%
\begin{equation}%
\begin{array}
[c]{ccccccccccc}%
A_{1} & = & -\dfrac{a-2b}{8c} &  & A_{2} & = & \dfrac{3(5a-6b)}{16c} &  &
A_{3} & = & -\dfrac{3(5a-6b)}{16c}%
\end{array}
\label{23a}%
\end{equation}
therefore%
\begin{equation}
\mathcal{S}=%
{\displaystyle\sum\limits_{n=0}^{\infty}}
\left(  \frac{27}{4}\rho\right)  ^{n}\left(  A_{1}\cdot\mathcal{B}_{n,1}%
+A_{2}\cdot\mathcal{B}_{n,2}+A_{3}\cdot\mathcal{B}_{n,3}\right) \label{23c}%
\end{equation}
where $\mathcal{B}_{n,k}$ is the Beta function integral
\begin{equation}
\mathcal{B}_{n,k}=B\left(  n+\tfrac{1}{2},2n+k\right)  =\frac{\Gamma
(2n+k)\,\Gamma\left(  n+\frac{1}{2}\right)  }{\Gamma\left(  3n+k+\frac{1}%
{2}\right)  }=%
{\displaystyle\int\limits_{0}^{1}}
(1-x)^{n-\frac{1}{2}}\,x^{2n+k-1}\,dx\label{23d}%
\end{equation}
by dominated convergence sum and integral in Eq.(\ref{23c}) commute, so that%
\begin{equation}
\mathcal{S=}%
{\displaystyle\int\limits_{0}^{1}}
\left(  A_{1}+A_{2}\,x+A_{3}\,x^{2}\right)  \cdot%
{\displaystyle\sum\limits_{n=0}^{\infty}}
\left(  \frac{27}{4}\rho\cdot x^{2}(1-x)\right)  ^{n}\frac{dx}{\sqrt{1-x}%
}\label{23e}%
\end{equation}
finally,%
\begin{equation}
\mathcal{S=}%
{\displaystyle\int\limits_{0}^{1}}
\frac{A_{1}+A_{2}\,x+A_{3}\,x^{2}}{1-\frac{27}{4}\rho\cdot x^{2}(1-x)}%
\cdot\frac{dx}{\sqrt{1-x}}\label{23f}%
\end{equation}
the proof follows through the evaluation of this integral. For integer
coefficients polynomials $u(x)$ and $v(x)$, this expression can be written as%
\begin{equation}
\mathcal{S=}%
{\displaystyle\int\limits_{0}^{1}}
\,\frac{u(x)}{v(x)}\cdot\frac{dx}{\sqrt{1-x}}\label{23g}%
\end{equation}
applying the $A_{1},A_{2},A_{3},\rho$ values included in \textit{Table} I
these polynomials are obtained \vspace*{4pt}

\begin{center}%
\begin{tabular}
[c]{||c|c|c||}\hline\hline
\textit{p} & $u(x)$ $%
\genfrac{}{}{0pt}{0}{\genfrac{}{}{0pt}{1}{{}}{{}}}{{}}%
$ & $%
\genfrac{}{}{0pt}{0}{\genfrac{}{}{0pt}{1}{{}}{{}}}{{}}%
v(x)$\\\hline\hline
{\small 2} & \multicolumn{1}{|r|}{$3x^{2}-3x+200%
\genfrac{}{}{0pt}{0}{\genfrac{}{}{0pt}{1}{{}}{{}}}{{}}%
$} & \multicolumn{1}{|r||}{$%
\genfrac{}{}{0pt}{0}{\genfrac{}{}{0pt}{1}{{}}{{}}}{{}}%
(x+8)(x^{2}-9x+72)$}\\
{\small 3} & \multicolumn{1}{|r|}{$2x^{2}-2x+20%
\genfrac{}{}{0pt}{0}{\genfrac{}{}{0pt}{1}{{}}{{}}}{{}}%
$} & \multicolumn{1}{|r||}{$%
\genfrac{}{}{0pt}{0}{\genfrac{}{}{0pt}{1}{{}}{{}}}{{}}%
(x+3)(x^{2}-4x+12)$}\\
{\small 5} & \multicolumn{1}{|r|}{$4x+16%
\genfrac{}{}{0pt}{0}{\genfrac{}{}{0pt}{1}{{}}{{}}}{{}}%
$} & \multicolumn{1}{|r||}{$%
\genfrac{}{}{0pt}{0}{\genfrac{}{}{0pt}{1}{{}}{{}}}{{}}%
x^{2}+4x+20$}\\
{\small 7} & \multicolumn{1}{|r|}{$972\,x^{2}-972\,x+840%
\genfrac{}{}{0pt}{0}{\genfrac{}{}{0pt}{1}{{}}{{}}}{{}}%
$} & \multicolumn{1}{|r||}{$%
\genfrac{}{}{0pt}{0}{\genfrac{}{}{0pt}{1}{{}}{{}}}{{}}%
(9x+7)(81x^{2}-144x+112)$}\\
{\small 10} & \multicolumn{1}{|r|}{$27x+45%
\genfrac{}{}{0pt}{0}{\genfrac{}{}{0pt}{1}{{}}{{}}}{{}}%
$} & \multicolumn{1}{|r||}{$%
\genfrac{}{}{0pt}{0}{\genfrac{}{}{0pt}{1}{{}}{{}}}{{}}%
9x^{2}+15x+40$}\\\hline\hline
\end{tabular}
\medskip
\end{center}

By using some CAS code to evaluate these integrals symbolically, the proof
terminates with these lines in Maple\texttrademark\ \cite{MAPLE}\vspace*{2pt}%
\begin{center}
\ifcase\msipdfoutput
\includegraphics[viewport=0in 0in 11.552200in 4.052500in,
height=1.9726in,
width=5.5728in
]%
{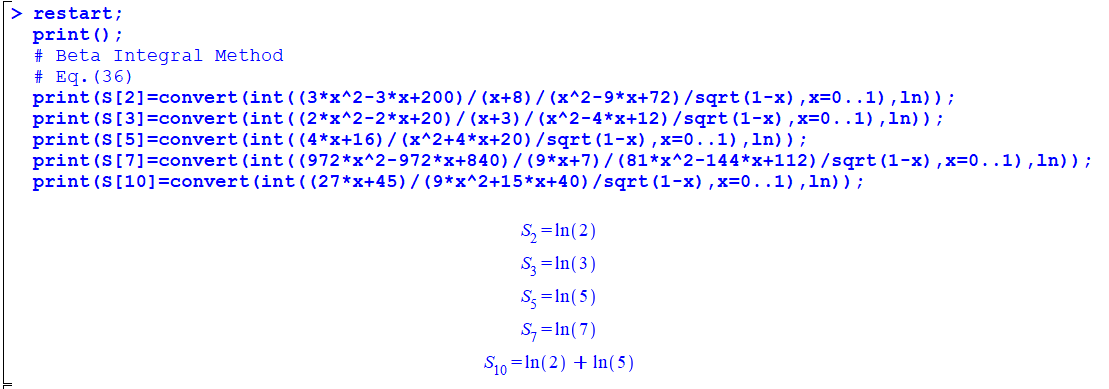}%
\else
\includegraphics[viewport=0in 0in 5.572800in 1.972600in,
height=1.9726in,
width=5.5728in
]%
{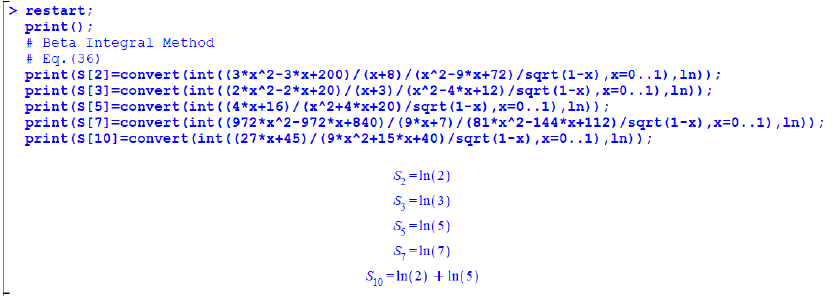}%
\fi
\end{center}

\subsection{Hypergeometric Closed Forms}

A different proof is derived from a recent work by Zhi-Wei Sun and Yajun Zhou
\cite{SUNZHOU} by means of some proven series closed forms.\ From Table 3,
rows A.1, B.1, C.1 and D.1 in this paper, these expressions are established%

\begin{equation}
\Phi_{A}(z)=%
{\displaystyle\sum\limits_{n=1}^{\infty}}
\frac{\binom{2n}{n}}{n\,\binom{3n}{n}\binom{6n}{3n}}\cdot\left[  \frac
{4}{z(1-z^{2})}\right]  ^{2n}=%
{\displaystyle\sum\limits_{n=1}^{\infty}}
\frac{\rho^{n}}{n}\cdot%
\begin{bmatrix}%
\genfrac{}{}{0pt}{0}{{}}{{}}%
1 & \frac{1}{2}%
\genfrac{}{}{0pt}{0}{{}}{{}}%
\\%
\genfrac{}{}{0pt}{0}{{}}{{}}%
\frac{1}{6} & \frac{5}{6}%
\genfrac{}{}{0pt}{0}{{}}{{}}%
\end{bmatrix}
_{n}\label{15}%
\end{equation}%
\[
\Phi_{B}(z)=%
{\displaystyle\sum\limits_{n=0}^{\infty}}
\frac{\binom{2n}{n}}{(6n+1)\,\binom{3n}{n}\binom{6n}{3n}}\cdot\left[  \frac
{4}{z(1-z^{2})}\right]  ^{2n}=%
{\displaystyle\sum\limits_{n=0}^{\infty}}
\frac{\rho^{n}}{6n+1}\cdot%
\begin{bmatrix}%
\genfrac{}{}{0pt}{0}{{}}{{}}%
1 & \frac{1}{2}%
\genfrac{}{}{0pt}{0}{{}}{{}}%
\\%
\genfrac{}{}{0pt}{0}{{}}{{}}%
\frac{1}{6} & \frac{5}{6}%
\genfrac{}{}{0pt}{0}{{}}{{}}%
\end{bmatrix}
_{n}%
\]%
\[
\Phi_{C}(z)=%
{\displaystyle\sum\limits_{n=0}^{\infty}}
\frac{\binom{2n}{n}}{(6n+5)\,\binom{3n}{n}\binom{6n}{3n}}\cdot\left[  \frac
{4}{z(1-z^{2})}\right]  ^{2n}=%
{\displaystyle\sum\limits_{n=0}^{\infty}}
\frac{\rho^{n}}{6n+5}\cdot%
\begin{bmatrix}%
\genfrac{}{}{0pt}{0}{{}}{{}}%
1 & \frac{1}{2}%
\genfrac{}{}{0pt}{0}{{}}{{}}%
\\%
\genfrac{}{}{0pt}{0}{{}}{{}}%
\frac{1}{6} & \frac{5}{6}%
\genfrac{}{}{0pt}{0}{{}}{{}}%
\end{bmatrix}
_{n}%
\]%
\[
\Phi_{D}(z)=%
{\displaystyle\sum\limits_{n=1}^{\infty}}
\frac{\binom{2n}{n}}{(2n-1)\,\binom{3n}{n}\binom{6n}{3n}}\cdot\left[  \frac
{4}{z(1-z^{2})}\right]  ^{2n}=%
{\displaystyle\sum\limits_{n=1}^{\infty}}
\frac{\rho^{n}}{2n-1}\cdot%
\begin{bmatrix}%
\genfrac{}{}{0pt}{0}{{}}{{}}%
1 & \frac{1}{2}%
\genfrac{}{}{0pt}{0}{{}}{{}}%
\\%
\genfrac{}{}{0pt}{0}{{}}{{}}%
\frac{1}{6} & \frac{5}{6}%
\genfrac{}{}{0pt}{0}{{}}{{}}%
\end{bmatrix}
_{n}%
\]
with%
\begin{equation}
\Phi_{A}(z)=\frac{1}{1-3z^{2}}\,\left[  \,6z\tanh^{-1}\left(  \frac{1}%
{z}\right)  +(3z^{2}-2)\cdot\Phi_{L}(z)\,\right] \label{16d}%
\end{equation}%
\[
\Phi_{B}(z)=\frac{3z(1-z^{2})}{2(1-3z^{2})}\,\left[  \tanh^{-1}\left(
\frac{1}{z}\right)  -\frac{z}{2}\cdot\Phi_{L}(z)\,\right]
\]%
\[
\Phi_{C}(z)=\frac{3z(1-z^{2})}{2(3z^{2}-1)}\left[  \,(9z^{4}-9z^{2}%
+1)\tanh^{-1}\left(  \frac{1}{z}\right)  +\frac{z\,(9z^{4}-15z^{2}+5)}%
{2}\,\cdot\Phi_{L}(z)\,\right]
\]%
\[
\Phi_{D}(z)=\frac{1}{2(1-z^{2})(1-3z^{2})}\left[  \,\frac{3z^{4}-3z^{2}+2}%
{z}\tanh^{-1}\left(  \frac{1}{z}\right)  +\frac{z^{2}(3z^{2}-5)}{2}\cdot
\Phi_{L}(z)\,\right]
\]

\bigskip

and%
\begin{equation}
\Phi_{L}(z)=i\cdot\frac{1}{\sqrt{3z^{2}-4}}\log\left(  \frac{z^{2}%
-2+i\sqrt{3z^{2}-4}}{z^{2}-2-i\sqrt{3z^{2}-4}}\right)  \vspace*{6pt}\label{16}%
\end{equation}
the series $\mathcal{S}$ are taken from the \textrm{RHS} of Eqs.(\ref{20}%
--\ref{21}), by setting
\begin{equation}
\rho=\frac{4}{27}\cdot\frac{1}{z^{2}(1-z^{2})^{2}}\label{16h}%
\end{equation}
For $\rho>0,$ $z$ is computed \vspace*{-6pt}as the real root of $z(1-z^{2}%
)=\tfrac{2}{3}\sqrt{\tfrac{1}{3\rho}}.$ For $\rho<0,$ $z=iy$ where $y$ is the
real root of $y(1+y^{2})=\tfrac{2}{3}\sqrt{\tfrac{-1}{3\rho}}.$ \textit{Table}
I, last column, has\vspace*{-2pt} the values of $z$ for each $\rho$. Linear
combinations of $\Phi_{A}$ and $\Phi_{D}$ or $\Phi_{B}$ and $\Phi_{C}$
generate $\mathcal{S}$ from either the \textrm{RHS} of Eq.(\ref{20}) or
Eq.(\ref{21}) respectively
\begin{equation}
\mathcal{S=-}\frac{\beta}{\gamma}\,\cdot\Phi_{A}(z)+\frac{\alpha+2\beta
}{\gamma}\cdot\Phi_{D}(z)\label{18a}%
\end{equation}%
\begin{equation}
\mathcal{S=}\ \frac{6b-a}{24c}\cdot\Phi_{B}(z)+\frac{5a-6b}{24c}\cdot\Phi
_{C}(z)\vspace*{8pt}\label{18b}%
\end{equation}
The following Maple\texttrademark\ code evaluates these expressions in
symbolic form and proves the identities\vspace*{-4pt}%

\begin{center}
\ifcase\msipdfoutput
\includegraphics[viewport=0in 0in 11.604000in 8.312600in,
height=4.849in,
width=6.7585in
]%
{./Log3_.png}%
\else
\includegraphics[viewport=0in 0in 6.758500in 4.849000in,
height=4.849in,
width=6.7585in
]%
{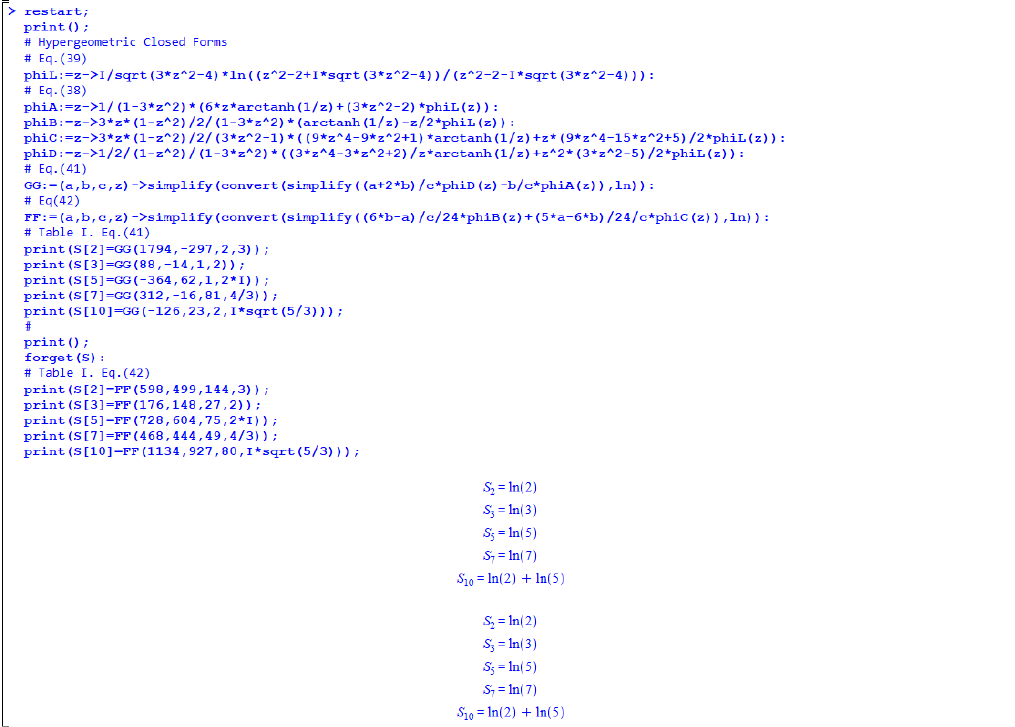}%
\fi
\end{center}
\vspace*{-10pt}

\subsection{Level 1 Ramanujan--type Series for Logarithms}

C. Hakimoglu's work \cite{CHBROWN} based on the Beta Integral method, can be
extended to get a general Ramanujan-type formula of level 1 (denominators of
signature 6) for logarithm convergent series.\vspace{-3pt} Starting there with
the proven equation 2.9 of Section 2.1, that converges for $|v|<\sqrt
{\tfrac{3}{2}}$ as\vspace*{-4pt}
\begin{equation}
-3\sqrt{1-v^{2}}\cdot\frac{\sin^{-1}v}{v}=%
{\displaystyle\sum\limits_{n=0}^{\infty}}
\left(  \frac{4\,v^{6}}{27\,(v^{2}-1)}\right)  ^{n}\cdot\frac{P_{0}%
(n,v)}{(6n+1)(6n+5)}\cdot%
\begin{bmatrix}%
\genfrac{}{}{0pt}{0}{{}}{{}}%
1 & \frac{1}{2}%
\genfrac{}{}{0pt}{0}{{}}{{}}%
\\%
\genfrac{}{}{0pt}{0}{{}}{{}}%
\frac{1}{6} & \frac{5}{6}%
\genfrac{}{}{0pt}{0}{{}}{{}}%
\end{bmatrix}
_{n}\vspace*{-4pt}\label{18w}%
\end{equation}%
\[
P_{0}(n,v)=(4v^{4}+6v^{2}-18)\cdot n+2v^{4}+5v^{2}-15\vspace*{0pt}%
\]
\ By setting $v=i\,u$, this identity is transformed using\vspace*{-2pt}%
\[
\sqrt{1-v^{2}}\cdot\frac{\sin^{-1}v}{v}=\frac{\sqrt{1+u^{2}}}{u}\cdot
\log\left(  u+\sqrt{1+u^{2}}\right)  \vspace*{-12pt}%
\]
so,\vspace*{-7pt}%
\begin{equation}
\log\left(  u+\sqrt{1+u^{2}}\right)  =-\frac{u}{3\sqrt{1+u^{2}}}\hspace*{2pt}%
{\displaystyle\sum\limits_{n=0}^{\infty}}
\left(  \frac{4\,u^{6}}{27\,(u^{2}+1)}\right)  ^{n}\hspace*{-4pt}\frac
{P_{1}(n,u)}{(6n+1)(6n+5)}%
\begin{bmatrix}%
\genfrac{}{}{0pt}{0}{{}}{{}}%
1 & \frac{1}{2}%
\genfrac{}{}{0pt}{0}{{}}{{}}%
\\%
\genfrac{}{}{0pt}{0}{{}}{{}}%
\frac{1}{6} & \frac{5}{6}%
\genfrac{}{}{0pt}{0}{{}}{{}}%
\end{bmatrix}
_{n}\vspace*{0pt}\label{18u}%
\end{equation}
\vspace*{-6pt}%
\[
P_{1}(n,u)=(4u^{4}-6u^{2}-18)\cdot n+2u^{4}-5u^{2}-15\vspace*{0pt}%
\]
By selecting $p=\left(  u+\sqrt{1+u^{2}}\right)  ^{2}$ such that $u=\pm
\,\frac{p-1}{2\sqrt{p}}$, the following identity is obtained\vspace*{0in}%

\begin{equation}
\log p=-\frac{1}{12\,}\cdot%
{\displaystyle\sum\limits_{n=0}^{\infty}}
\frac{(p-1)^{6n+1}}{108^{n}\,p^{2n+2}\,(p+1)^{2n+1}}\cdot\frac{P(n,p)}%
{(6n+1)(6n+5)}\cdot%
\begin{bmatrix}%
\genfrac{}{}{0pt}{0}{{}}{{}}%
1 & \frac{1}{2}%
\genfrac{}{}{0pt}{0}{{}}{{}}%
\\%
\genfrac{}{}{0pt}{0}{{}}{{}}%
\frac{1}{6} & \frac{5}{6}%
\genfrac{}{}{0pt}{0}{{}}{{}}%
\end{bmatrix}
_{n}\vspace*{2pt}\label{18z}%
\end{equation}%
\[
P(n,p)=2\,(p^{2}-14p+1)(p^{2}+4p+1)\cdot n+p^{4}-14p^{3}-94p^{2}%
-14p+1\vspace*{6pt}%
\]
being valid for $|p-7|<4\sqrt{3}.$ This means that for $p\in%
\mathbb{Z}
\cap\lbrack2,13]$ a monotone linearly convergent logarithm series with
$\rho\in\lbrack0,1]$ is obtained beyond the cases $p=2,3,7$ displayed in
\textit{Table} I. Note that this proves logarithm series that have a
convergence rate%
\begin{equation}
\rho=\frac{(p-1)^{6}}{108\,p^{2}\,(p+1)^{2}}\label{18y}%
\end{equation}
producing the fast series Eqs.(\ref{8}--\ref{8a}) for the smaller values
$p=2,3$ and also $p=7$. Eq.(\ref{18z}) can be compared to Z.W. Sun's
\cite{ZHINEW} Corollary 1.5, where level 2 (denominators of signature 4)
Ramanujan type series valid for $p_{0}^{-1}$ $<p<$\ $p_{0}=$
{\small 21.26668...}, are derived as %
(Notice that a level 3 series -not shown- can be derived from this reference as well)
\begin{equation}
\log p=-\frac{1}{2\,}\cdot%
{\displaystyle\sum\limits_{n=0}^{\infty}}
\frac{(p-1)^{4n+1}}{16^{n}\,(-p)^{n+1}\,(p+1)^{2n+1}}\cdot\frac{Q(n,p)}%
{(4n+1)(4n+3)}\cdot%
\begin{bmatrix}%
\genfrac{}{}{0pt}{0}{{}}{{}}%
1 & \frac{1}{2}%
\genfrac{}{}{0pt}{0}{{}}{{}}%
\\%
\genfrac{}{}{0pt}{0}{{}}{{}}%
\frac{1}{4} & \frac{3}{4}%
\genfrac{}{}{0pt}{0}{{}}{{}}%
\end{bmatrix}
_{n}\label{18m}%
\end{equation}%
\[
Q(n,p)=2(p^{2}+6p+1)\cdot n+p^{2}+10p+1\vspace*{6pt}%
\]
Eq.(\ref{18z}) and Eq.(\ref{18m}) are also found proven by WZ pairs
Eq.(\ref{7}--\ref{8}) with transformations $\left.  W\hspace*{-2pt}%
Z_{2,1}^{\text{{\tiny (1)}}}\right\vert _{p}$\thinspace or\thinspace$\left.
W\hspace*{-2pt}Z_{1,2}^{\text{{\tiny (2)}}}\right\vert _{p}$\thinspace
and\thinspace$\left.  W\hspace*{-2pt}Z_{1,1}^{\text{{\tiny (1)}}}\right\vert
_{p}$\thinspace or\thinspace$\left.  W\hspace*{-2pt}Z_{1,1}^{\text{{\tiny (2)}%
}}\right\vert _{p}$\thinspace respectively. $\vspace*{-2pt}$ \hspace{-4pt} In addition, the next Maple\texttrademark code gives a proof of the generalized identity
Eq.(\ref{18z}) applying both methods Beta Integral (BI) and\ hypergeometric
Closed Forms (CF). Note that solving Eq.(\ref{16h}) for $z$ with $\rho$\ from
Eq.(\ref{18y}) produces $z=\pm\frac{p+1}{p-1}$ which gives the rightmost
column values of Table I if $\rho>0$,\vspace*{0pt}
\begin{center}
\ifcase\msipdfoutput
\includegraphics[viewport=0in 0in 11.854000in 7.062900in,
height=3.2768in,
width=5.4812in
]%
{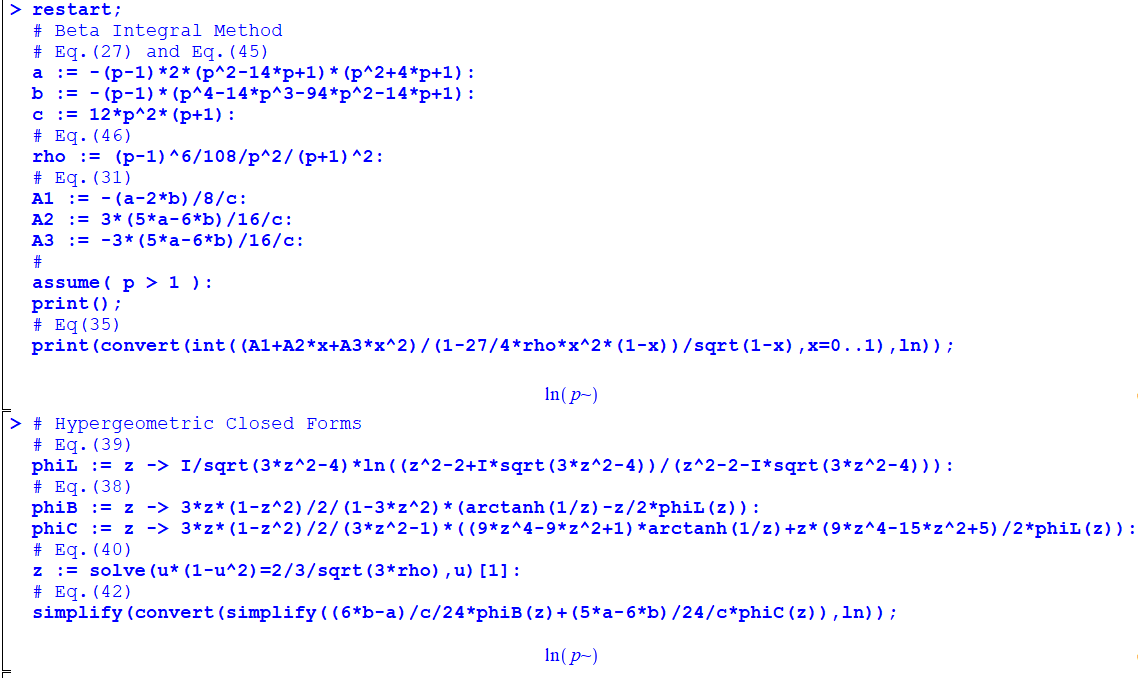}%
\else
\includegraphics[viewport=0in 0in 5.481200in 3.276800in,
height=3.2768in,
width=5.4812in
]%
{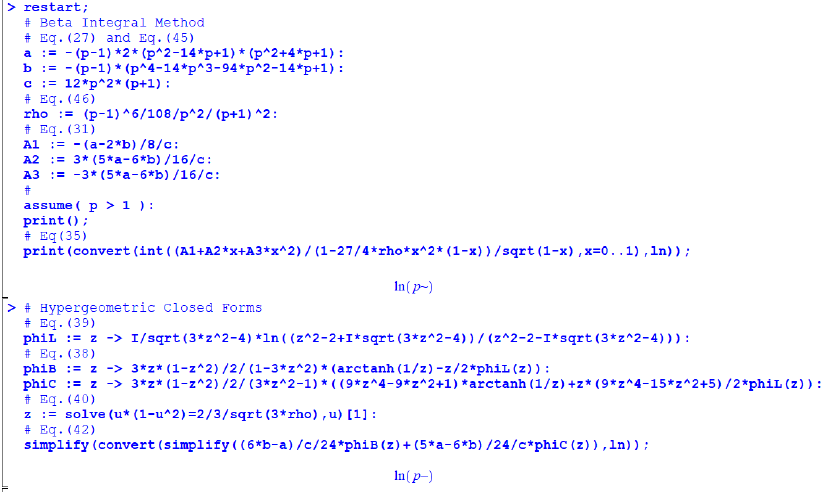}%
\fi
\end{center}
\vspace*{-6pt}\vspace*{4pt}%
\begin{center}
\ifcase\msipdfoutput
\includegraphics[viewport=0in 0in 3.333000in 3.333000in,
trim=-0.010426in 0.031277in 0.010426in -0.031277in,
height=3.3615in,
width=3.3615in
]%
{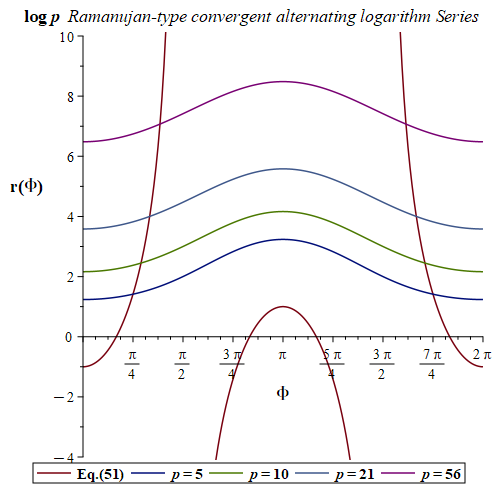}%
\else
\includegraphics[viewport=0in 0in 3.361500in 3.361500in,
height=3.3615in,
width=3.3615in
]%
{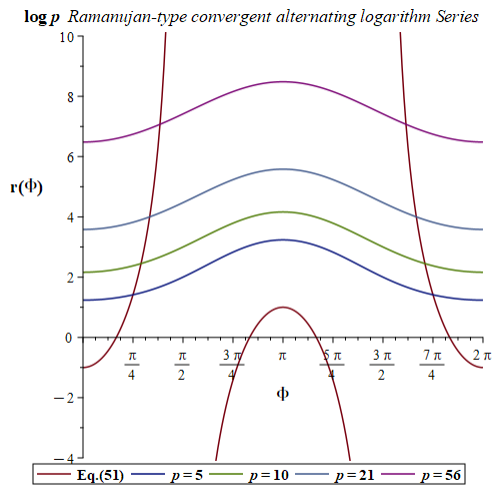}%
\fi
\end{center}
\vspace*{0in}

\subsubsection{Alternating Series}

For all real $p>0$ Eq.(\ref{18y}) produces only monotone series with $\rho
\geq0.$ To get alternating series like $\log10$ and the super fast $\log5$
displayed in \textit{Table} I it is necessary to work in the complex quadratic
number field $%
\mathbb{Q}
(\sqrt{m})$, $m\in%
\mathbb{Z}
_{<0},$ $m$ squarefree,\ forcing $\arg(\rho)=(2k+1)\pi$ with $k\in%
\mathbb{Z}
$ and $p$ expressed as the norm\vspace*{0pt}
\begin{equation}
p=(1+re^{i\phi})(1+re^{-i\phi})=1+2r\cos\phi+r^{2}\vspace*{-8pt}\label{19q}%
\end{equation}
so%
\begin{equation}
\log p=\log(1+re^{i\phi})+\log(1+re^{-i\phi})=2\operatorname{Re}\left(
\log(1+re^{i\phi})\right)  \vspace*{4pt}\label{19t}%
\end{equation}
The complex logarithm series $\log(1+re^{i\phi})$ is computed using identities
Eq.(\ref{18z}--\ref{18y}) but compelling the convergence rate\vspace*{-6pt}
\begin{equation}
\rho(r,\phi)=\frac{r^{6}e^{6i\phi}}{108\,(1+re^{i\phi})^{2}\,(2+re^{i\phi
})^{2}}\vspace*{-2pt}\label{19p}%
\end{equation}
to give alternating series for $\phi\in\lbrack0,2\pi)$. This means\vspace
*{-2pt}%
\[
3\phi-\arg(1+re^{i\phi})\,-\arg(2+re^{i\phi})=(2k+1)\frac{\pi}{2},\text{
\ \ }k\in%
\mathbb{Z}
\vspace*{-4pt}%
\]
that is equivalent to\vspace*{-4pt}%
\begin{equation}
r^{2}\cdot\cos\phi+3r\cdot\cos2\phi+2\cos3\phi=0\vspace*{4pt}\label{19a}%
\end{equation}
this is the red curve $r(\phi)$ in the graph above. This identity together
with,%
\begin{equation}
r^{2}+2r\cdot\cos\phi+1-p=0\label{19b}%
\end{equation}
give a system of two equations with input parameter $p$ and variables
$(r,\phi)$ to be solved. For fixed $p>1$ this system has only one solution
with $r>0$ in the range $\phi\in\left(  0,\pi\right)  $. To get the
convergence range for $p,$ $\rho(r,\phi)=-1$ --Eq.(\ref{19p})-- must be added,
this is\vspace*{2pt}%
\begin{equation}
r^{6}-108\cdot\left(  r^{2}+2r\cdot\cos\phi+1\right)  \cdot\left(
r^{2}+4r\cdot\cos\phi+4\right)  =0\label{20b}%
\end{equation}%
\[%
\begin{tabular}
[c]{ccccccccc}\hline\hline
\multicolumn{1}{||c}{$%
\genfrac{}{}{0pt}{0}{\genfrac{}{}{0pt}{1}{{}}{{}}}{{}}%
p$} & \multicolumn{1}{|c}{$m$} & \multicolumn{1}{|c}{$r$} &
\multicolumn{1}{|c}{$\phi$} & \multicolumn{1}{|c}{$1+re^{i\phi}$} &
\multicolumn{1}{|c}{$a$} & \multicolumn{1}{|c}{$b$} & \multicolumn{1}{|c}{$c$}
& \multicolumn{1}{|c||}{$\rho(r,\phi)$}\\\hline\hline
\multicolumn{1}{||c}{$%
\genfrac{}{}{0pt}{0}{\genfrac{}{}{0pt}{1}{{}}{{}}}{{}}%
${\small 5}} & \multicolumn{1}{|c}{$-1$} & \multicolumn{1}{|c}{$\sqrt{2}$} &
\multicolumn{1}{|c}{$\frac{\pi}{4}$} & \multicolumn{1}{|c}{{\small 2 +
}\textit{i}} & \multicolumn{1}{|r}{{\small 728}} &
\multicolumn{1}{|r}{{\small 604}} & \multicolumn{1}{|r}{{\small 75}} &
\multicolumn{1}{|c||}{{\small --}$\frac{1}{675}$}\\\hline
\multicolumn{1}{||c}{$%
\genfrac{}{}{0pt}{0}{\genfrac{}{}{0pt}{1}{{}}{{}}}{{}}%
${\small 10}} & \multicolumn{1}{|c}{$-15$} & \multicolumn{1}{|c}{$\sqrt{6}$} &
\multicolumn{1}{|c}{$\tan^{-1}\sqrt{\frac{5}{3}}$} &
\multicolumn{1}{|c}{$\frac{5+i\sqrt{15}}{2}$} &
\multicolumn{1}{|r}{{\small 1134}} & \multicolumn{1}{|r}{{\small 927}} &
\multicolumn{1}{|r}{{\small 80}} & \multicolumn{1}{|c||}{{\small --}$\frac
{1}{80}$}\\\hline
\multicolumn{1}{||c}{$%
\genfrac{}{}{0pt}{0}{\genfrac{}{}{0pt}{1}{{}}{{}}}{{}}%
${\small 21}} & \multicolumn{1}{|c}{$-3$} & \multicolumn{1}{|c}{$4$} &
\multicolumn{1}{|c}{$\frac{\pi}{3}$} & \multicolumn{1}{|c}{{\small 3 +
2}\textit{i}$\sqrt{{\small 3}}$} & \multicolumn{1}{|r}{{\small 8840}} &
\multicolumn{1}{|r}{{\small 6940}} & \multicolumn{1}{|r}{{\small 441}} &
\multicolumn{1}{|c||}{{\small --}$\frac{256}{3969}$}\\\hline
\multicolumn{1}{||c}{$%
\genfrac{}{}{0pt}{0}{\genfrac{}{}{0pt}{1}{{}}{{}}}{{}}%
${\small 56}} & \multicolumn{1}{|c}{$-7$} & \multicolumn{1}{|c}{$5\sqrt{2}$} &
\multicolumn{1}{|c}{$\tan^{-1}\sqrt{7}$} & \multicolumn{1}{|c}{$\frac
{7+5i\sqrt{7}}{2}$} & \multicolumn{1}{|r}{{\small 179630}} &
\multicolumn{1}{|r}{{\small 126775}} & \multicolumn{1}{|r}{{\small 5376}} &
\multicolumn{1}{|c||}{{\small --}$\frac{15625}{48384}$}\\\hline\hline
&  &  &  &  &  &  &  & \\
\multicolumn{9}{c}{$\vspace*{4pt}$\textit{Table}{\small \ }II. \textit{Eq.}%
{\small (\ref{21})} \textit{Ramanujan--type alternating logarithm series}}%
\end{tabular}
\vspace*{4pt}%
\]

Solving Eqs.(\ref{19a}--\ref{19b}--\ref{20b}) it gives $(r_{L},\phi_{L}%
,p_{L})=$ {\normalsize {\small (11.2691..., 1.3233..., 133.5126...), t}}he
convergence numerical limits, so that alternating series go up to integers
$p<134${\normalsize {\small \ }}larger than the monotone series convergence
integer limit\ $p<14.$ The whole range $p\in%
\mathbb{Z}
\cap\lbrack2,133]$ was scanned and Eqs.(\ref{19a}--\ref{19b}) were solved to
get $(r,\phi)$ for each $p$. Almost all integers in this range give solutions
leading up to irreducible algebraic $\rho(r,\phi)$ which are discarded since
they do not bring rational series. Only four sporadic alternating series are
good having rational expressions with solutions in $%
\mathbb{Q}
(\sqrt{m})$, $m\in%
\mathbb{Z}
_{<0},$ they are $p=5,10,21$ and $56$ where just $p=5$ is a fast series. This
identity from Eq.(\ref{21}), Eq.(\ref{18z}) and Eq.(\ref{19t}),
\begin{equation}
\frac{a}{c}\cdot n+\frac{b}{c}=-\frac{1}{6}\operatorname{Re}\left(  \left.
\frac{(p-1)}{p^{2}(p+1)}\cdot P(n,p)\right\vert _{p=1+re^{i\phi}}\right)
\label{20a}%
\end{equation}
is applied for these four cases to give the integers $a,b,c$. The graph
$r(\phi)$ in the previous page shows bell shaped curves, they are,
bottom--top, $(p,m)=(5,-1)${\small \textit{\ blue}}, $(10,-15)$%
{\small \textit{\ green}}, $(21,-3)${\small \textit{\ light-blue}},
$(56,-7)${\small \textit{\ purple.}} They cut the curve Eq.(\ref{19a}%
){\small \textit{\ red,}} on the left giving crosspoints $(r,\phi),$ the
quadratic field identifier $m $ and by Eq.(\ref{21}) and Eq.(\ref{20a}) the
integers $a,b,c$ that are displayed in \textit{Table }II above\vspace*{6pt}.

\subsection{$W\hspace*{-2pt}Z_{s,t}$ \ Certificates}

Except for Eq.(\ref{13}) all series in Section 5 have $W\hspace*{-2pt}Z_{s,t}
$ certificates brought by the WZ algorithm providing another proof of such
identities. For $d=2$ series Eqs.(\ref{8}--\ref{8a}--\ref{8b}) are
$W\hspace*{-2pt}Z_{s,t}$ proven. Applying Eqs.(\ref{4a}--\ref{7}),
certificates%
\begin{equation}
R_{s,t}(n,k)=\frac{G_{s,t}(n,k)}{F_{s,t}(n,k)}.\vspace*{-6pt}\label{18c}%
\end{equation}
are, for Eq.(\ref{8}),%
\begin{equation}%
\begin{array}
[c]{ccc}%
R_{2,1}^{\text{{\tiny (1)}}}(n,k)|_{p=2} & = & \dfrac{144k^{2}%
+(828n+558)k+1196n^{2}+1596n+499}{32(6n+2k+3)(6n+5+2k)}\\
&  & \\
R_{1,2}^{\text{{\tiny (2)}}}(n,k)|_{p=2} & = & \dfrac{128k^{2}%
+(782n+519)k+1196n^{2}+1596n+499}{36(6n+2k+3)(6n+5+2k)}%
\end{array}
\label{18d}%
\end{equation}
for Eq.(\ref{8a}),%
\begin{equation}%
\begin{array}
[c]{ccc}%
R_{2,1}^{\text{{\tiny (1)}}}(n,k)|_{p=3} & = & \dfrac{48k^{2}%
+(256n+176)k+352n^{2}+472n+148}{9(6n+2k+3)(6n+5+2k)}\\
&  & \\
R_{1,2}^{\text{{\tiny (2)}}}(n,k)|_{p=3} & = & \dfrac{9k^{2}+(56n+37)k+88n^{2}%
+118n+37}{3(6n+2k+3)(6n+5+2k)}%
\end{array}
\vspace*{4pt}\label{18e}%
\end{equation}
and for Eq.(\ref{8b}),
\begin{equation}%
\begin{array}
[c]{ccl}%
R_{2,1}^{\text{{\tiny (1)}}}(n,k)|_{p=2\pm i} & =%
\genfrac{}{}{0pt}{0}{\genfrac{}{}{0pt}{1}{{}}{{}}}{{}}%
& \mp i\cdot\dfrac{\left(  10k+26n+23\right)  \left(  2k+2n+1\right)
}{25(6n+2k+3)(6n+5+2k)}+\vspace*{6pt}\\
&
\genfrac{}{}{0pt}{0}{\genfrac{}{}{0pt}{1}{{}}{{}}}{{}}%
& +\dfrac{110k^{2}+(646n+433)k+936n^{2}+1246n+389}{25(6n+2k+3)(6n+5+2k)}%
\vspace*{6pt}\\
R_{1,2}^{\text{{\tiny (2)}}}(n,k)|_{p=2\pm i} & =%
\genfrac{}{}{0pt}{0}{\genfrac{}{}{0pt}{1}{{}}{{}}}{{}}%
& \mp i\cdot\dfrac{2(8k+26n+21)(k+2n+1)}{25(6n+5+2k)(6n+2k+3)}+\vspace*{6pt}\\
&
\genfrac{}{}{0pt}{0}{\genfrac{}{}{0pt}{1}{{}}{{}}}{{}}%
& +\dfrac{88k^{2}+(542n+359)k+832n^{2}+1108n+346}{25(6n+5+2k)(6n+2k+3)}%
\end{array}
\vspace*{4pt}\label{18f}%
\end{equation}
Other $d=4,6$ series have also $W\hspace*{-2pt}Z_{s,t}$ proofs but there is no
room here to put their $R_{s,t}$ certificates. In these cases, however, an
extension of the Beta Integral method can be applied to get simpler human
proofs as it is seen in the following\vspace*{0pt}.

\subsection{Beta Function Method Extended}

For particular hypergeometric motives some series can be transformed into
integrals having polynomials ratio integrands with simple algebraic kernels.
This is done splitting the series summands by partial fractions and applying
the Beta function integral form. By using symbolic integration packages, for
example either Mathematica Rubi \cite{RUBI} or Maple\texttrademark, these
integrals can be evaluated in closed form and the hypergeometric series are
demonstrated. From Eq.(\ref{2}) some identities for $d\geq2$ are proven in
this way using a similar method to the one applied on Eqs.(\ref{22}%
\thinspace--\thinspace\ref{23g})\vspace*{4pt}%
\begin{equation}%
\begin{tabular}
[c]{ccccc}%
$\mathcal{M}(n)$ & $=$ & $\left[
\begin{array}
[c]{cccc}%
r_{1} & r_{2} & ... & r_{d}\\
q_{1} & q_{2} & ... & q_{d}%
\end{array}
\vspace*{1pt}\right]  _{n}$ & $=$ & $4^{\nu n}\cdot\lambda^{n}\cdot\dfrac
{(\nu\hspace*{1pt}n)!\,(2m\hspace*{1pt}n)!\,(N\hspace*{1pt}n)!}{(m\hspace
*{1pt}n)!\,(2N\hspace*{1pt}n)!}$\\
&  &  &  & \\
& $=$ & $4^{\nu n}\cdot\lambda^{n}\cdot\dfrac{\dbinom{2mn}{mn}}{\dbinom
{Nn}{mn}\dbinom{2Nn}{Nn}}$ & $=$ & $\lambda^{n}\cdot\dfrac{\Gamma(\nu
\hspace*{1pt}n+1)\,\Gamma(m\hspace*{1pt}n+\frac{1}{2})}{\Gamma(N\hspace
*{1pt}n+\frac{1}{2})}$%
\end{tabular}
\tag{59}\label{59}%
\end{equation}
where $m\in%
\mathbb{Z}
_{\geq0},\ \nu\in%
\mathbb{Z}
_{>0},$ $N=\nu+m$ and $\lambda=\frac{N^{N}}{m^{m}\nu^{\nu}}$ with convention
$0^{0}=1.$ This gives the following table\vspace*{3pt}
\[%
\begin{tabular}
[c]{|c|c|c|c|c|c|c|}\hline\hline
$m$ & $\nu$ & $N$ & $\lambda$ & $\left[
\begin{array}
[c]{cccc}%
r_{1} & r_{2} & ... & r_{d}%
\end{array}
\right]  $ & $\left[
\begin{array}
[c]{cccc}%
q_{1} & q_{2} & ... & q_{d}%
\end{array}
\right]  $ & $%
\genfrac{}{}{0pt}{0}{\genfrac{}{}{0pt}{1}{{}}{{}}}{{}}%
\genfrac{}{}{0pt}{}{\text{Proven}}{\text{Identity}}%
\genfrac{}{}{0pt}{0}{\genfrac{}{}{0pt}{1}{{}}{{}}}{{}}%
$\\\hline\hline
{\small 1} & {\small 1} & {\small 2} & {\small 4} & $\left[  1\ \frac{1}%
{2}\right]  $ & $\left[  \frac{1}{4}\ \frac{3}{4}\right]  $ & $%
\genfrac{}{}{0pt}{0}{\genfrac{}{}{0pt}{1}{{}}{{}}}{{}}%
${\small Eq.(\ref{18m})}$%
\genfrac{}{}{0pt}{0}{\genfrac{}{}{0pt}{1}{{}}{{}}}{{}}%
$\\\hline
{\small 1} & {\small 2} & {\small 3} & $\frac{27}{4}$ & $\left[  1\ \frac
{1}{2}\right]  $ & $\left[  \frac{1}{6}\ \frac{5}{6}\right]  $ & $%
\genfrac{}{}{0pt}{0}{\genfrac{}{}{0pt}{1}{{}}{{}}}{{}}%
${\small Eq.(\ref{18z})}$%
\genfrac{}{}{0pt}{0}{\genfrac{}{}{0pt}{1}{{}}{{}}}{{}}%
$\\\hline
{\small 1} & {\small 4} & {\small 5} & $\frac{3125}{256}$ & $\left[
1\ \frac{1}{2}\ \frac{1}{4}\ \frac{3}{4}\right]  $ & $\left[  \frac{1}%
{10}\ \frac{3}{10}\ \frac{7}{10}\ \frac{9}{10}\right]  $ & $%
\genfrac{}{}{0pt}{0}{\genfrac{}{}{0pt}{1}{{}}{{}}}{{}}%
${\small Eq.(\ref{11})}$%
\genfrac{}{}{0pt}{0}{\genfrac{}{}{0pt}{1}{{}}{{}}}{{}}%
$\\\hline
{\small 3} & {\small 2} & {\small 5} & $\frac{3125}{108}$ & $\left[
1\ \frac{1}{2}\ \frac{1}{6}\ \frac{5}{6}\right]  $ & $\left[  \frac{1}%
{10}\ \frac{3}{10}\ \frac{7}{10}\ \frac{9}{10}\right]  $ & $%
\genfrac{}{}{0pt}{0}{\genfrac{}{}{0pt}{1}{{}}{{}}}{{}}%
${\small Eq.(\ref{9}), Eq.(\ref{15a})}$%
\genfrac{}{}{0pt}{0}{\genfrac{}{}{0pt}{1}{{}}{{}}}{{}}%
$\\\hline
{\small 3} & {\small 4} & {\small 7} & $\frac{823543}{6912}$ & $\left[
1\ \frac{1}{2}\ \frac{1}{4}\ \frac{3}{4}\ \frac{1}{6}\ \frac{5}{6}\right]  $ &
$\left[  \frac{1}{14}\ \frac{3}{14}\ \frac{5}{14}\ \frac{9}{14}\ \frac{11}%
{14}\ \frac{13}{14}\right]  $ & $%
\genfrac{}{}{0pt}{0}{\genfrac{}{}{0pt}{1}{{}}{{}}}{{}}%
${\small Eq.(\ref{18})}$%
\genfrac{}{}{0pt}{0}{\genfrac{}{}{0pt}{1}{{}}{{}}}{{}}%
$\\\hline\hline
\end{tabular}
\vspace*{6pt}%
\]
Note that these $(m,\nu)$ BI demonstrations correspond to $W\hspace
*{-2pt}Z_{\nu,m}^{\text{{\tiny (1)}}}$ or $W\hspace*{-2pt}Z_{m,\nu
}^{\text{{\tiny (2)}}}$ transformations applied on Eqs.(\ref{6}--\ref{7}),
namely these BI and WZ proofs are equivalent on the same series. For $d=4$ the
BI proofs details are left as an exercise to the interested reader. The next
example proving Eq.(\ref{18}) with{\small \ }$d=6$ gives the corresponding guidelines.

\subsubsection{Example I}

BI proof of Eq.(\ref{18}). Series is first converted by shifting the sum index
to $n=0.$ This gives \vspace*{-12pt}%

\begin{equation}%
\begin{tabular}
[c]{rrr}%
$\log2$ & $=$ & $%
{\displaystyle\sum\limits_{n=0}^{\infty}}
\dfrac{P(n)}{R(n)}\cdot\left(  \dfrac{1}{2^{4}\cdot3^{3}\cdot7^{7}}\right)
^{n}\cdot%
\begin{bmatrix}%
\genfrac{}{}{0pt}{0}{{}}{{}}%
1 & \frac{1}{2} & \frac{1}{4} & \frac{3}{4} & \frac{1}{6} & \frac{5}{6}%
\genfrac{}{}{0pt}{0}{{}}{{}}%
\\%
\genfrac{}{}{0pt}{0}{{}}{{}}%
\frac{1}{14} & \frac{3}{14} & \frac{5}{14} & \frac{9}{14} & \frac{11}{14} &
\frac{13}{14}%
\genfrac{}{}{0pt}{0}{{}}{{}}%
\end{bmatrix}
_{n}$%
\end{tabular}
\tag{60}\label{60}%
\end{equation}
where polynomials $P(n)$ and $R(n)$ are now%
\[%
\begin{tabular}
[c]{lll}%
{\small \textit{P}(\textit{n})} & $=$ & {\small 81969540480\thinspace
\textit{n}}$^{5}$ {\small + 239897521920\thinspace\textit{n}}$^{4}$ {\small +
266389817304\thinspace\textit{n}}$^{3}%
\genfrac{}{}{0pt}{0}{\genfrac{}{}{0pt}{1}{{}}{{}}}{{}}%
$\\
&  & {\small + 138594927588\thinspace\textit{n}}$^{2}$ {\small +
33273401586\thinspace\textit{n }+ 2913463287}$%
\genfrac{}{}{0pt}{0}{\genfrac{}{}{0pt}{1}{{}}{{}}}{{}}%
$\\
{\small \textit{R}(\textit{n})} & $=$ & {\small 217728\thinspace
(14\textit{n\thinspace}+\textit{\thinspace}1)(14\textit{n\thinspace
}+\textit{\thinspace}3)(14\textit{n\thinspace}+\textit{\thinspace
}5)(14\textit{n\thinspace}+\textit{\thinspace}9)(14\textit{n\thinspace
}+\textit{\thinspace}11)(14\textit{n\thinspace}+\textit{\thinspace}13)}$%
\genfrac{}{}{0pt}{0}{\genfrac{}{}{0pt}{1}{{}}{{}}}{{}}%
$%
\end{tabular}
\vspace*{0pt}%
\]

Take $\mathcal{G}%
(n)=P(n)/R(n)$ and $m=3$,$\ \nu=4$,$\ N=7$ from previous page table, %
 partial fractions PARI GP custom script \texttt{pfbeta(}%
$\mathcal{G}$\texttt{,7,}$\nu$\texttt{,1,}$N$\texttt{,1/2,0)} provides the following
coefficients (see help
\texttt{%
	$>$%
	?pfbeta}, code is available as it is indicated in Section 8.)\vspace*{2pt}
\[%
\begin{tabular}
[c]{lllllll}%
$A_{1}=\frac{3}{8}$ & $%
\genfrac{}{}{0pt}{0}{\genfrac{}{}{0pt}{1}{{}}{{}}}{{}}%
$ & $A_{2}=-\frac{563}{12096}$ & $%
\genfrac{}{}{0pt}{0}{\genfrac{}{}{0pt}{1}{{}}{{}}}{{}}%
$ & $A_{3}=\frac{479}{96768}$ & $%
\genfrac{}{}{0pt}{0}{\genfrac{}{}{0pt}{1}{{}}{{}}}{{}}%
$ & $A_{4}=-\frac{17}{110592}$\\
$A_{5}=\frac{91}{995328}$ & $%
\genfrac{}{}{0pt}{0}{\genfrac{}{}{0pt}{1}{{}}{{}}}{{}}%
$ & $A_{6}=-\frac{11}{995328}$ & $%
\genfrac{}{}{0pt}{0}{\genfrac{}{}{0pt}{1}{{}}{{}}}{{}}%
$ & $A_{7}=\frac{1}{995328}$ & $%
\genfrac{}{}{0pt}{0}{\genfrac{}{}{0pt}{1}{{}}{{}}}{{}}%
$ &
\end{tabular}
\]
such that the summands rational part is splitted as\vspace*{-3pt}%
\begin{equation}
\mathcal{G}(n)\text{ \ }=\text{\ \ }\frac{P(n)}{R(n)}\text{ \ }=\text{\ \ }%
\sum_{k=1}^{N}\frac{(\nu\,n+1)_{k-1}}{(N\,n+\frac{1}{2})_{k}}\cdot
A_{k}\tag{61}\label{60a}%
\end{equation}
Therefore using Eq.(\ref{59}) and $\rho=2^{-4}3^{-3}7^{-7}$ the RHS sum
$\mathcal{S}$ in Eq.(\ref{60}) can be written as%
\[%
\begin{tabular}
[c]{lcl}%
$\mathcal{S}$ & $=$ & $%
{\displaystyle\sum\limits_{n=0}^{\infty}}
\rho^{n}\cdot\lambda^{n}\cdot\dfrac{\Gamma(m\hspace*{1pt}n+\frac{1}%
{2})\,\Gamma(\nu\hspace*{1pt}n+1)\,}{\Gamma(N\hspace*{1pt}n+\frac{1}{2})}\cdot%
{\displaystyle\sum\limits_{k=1}^{N}}
\dfrac{(\nu\,n+1)_{k-1}}{(N\,n+\frac{1}{2})_{k}}\cdot A_{k}%
\genfrac{}{}{0pt}{0}{\genfrac{}{}{0pt}{1}{{}}{{}}}{{}}%
$\\
& $\mathcal{=}$ & $%
{\displaystyle\sum\limits_{n=0}^{\infty}}
\rho^{n}\cdot\lambda^{n}\cdot%
{\displaystyle\sum\limits_{k=1}^{N}}
A_{k}\cdot B(m\hspace*{1pt}n+\frac{1}{2},\nu\hspace*{1pt}n+k)%
\genfrac{}{}{0pt}{0}{\genfrac{}{}{0pt}{1}{{}}{{}}}{{}}%
$\\
& $\mathcal{=}$ & $%
{\displaystyle\sum\limits_{n=0}^{\infty}}
\rho^{n}\cdot\lambda^{n}\cdot%
{\displaystyle\sum\limits_{k=1}^{N}}
A_{k}\cdot%
{\displaystyle\int\nolimits_{0}^{1}}
(1-x)^{mn-1/2}\,x^{\nu n+k-1}dx%
\genfrac{}{}{0pt}{0}{\genfrac{}{}{0pt}{1}{{}}{{}}}{{}}%
$%
\end{tabular}
\]
using the Beta function. By dominated convergence sum and integral commute
giving
\begin{equation}%
\begin{tabular}
[c]{l}%
$\mathcal{S}\text{ \ }\mathcal{=}\text{ \ }%
{\displaystyle\int\nolimits_{0}^{1}}
\dfrac{\mathcal{P}_{N-1}(x)}{1-\rho\cdot\lambda\cdot x^{\nu}(1-x)^{m}}%
\cdot\dfrac{dx}{\sqrt{1-x}}%
\genfrac{}{}{0pt}{0}{\genfrac{}{}{0pt}{1}{{}}{{}}}{{}}%
$\\
$\mathcal{P}_{N-1}(x)\text{ \ }\mathcal{=}\text{ \ }A_{1}+A_{2}\,x+...+A_{N}%
\,x^{N-1}%
\genfrac{}{}{0pt}{0}{\genfrac{}{}{0pt}{1}{\genfrac{}{}{0pt}{1}{{}}{{}}}{{}%
}}{{}}%
$%
\end{tabular}
\tag{62}\label{60b}%
\end{equation}
Finally the following Maple lines complete the proof,\vspace*{-3pt}%
\begin{center}
\ifcase\msipdfoutput
\includegraphics[viewport=0in 0in 8.937800in 1.906100in,
height=1.2445in,
width=5.7692in
]%
{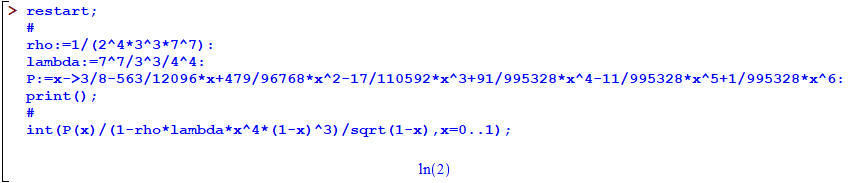}%
\else
\includegraphics[viewport=0in 0in 5.769200in 1.244500in,
height=1.2445in,
width=5.7692in
]%
{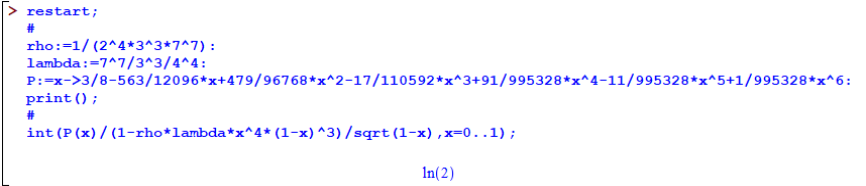}%
\fi
\end{center}

\section{log$\hspace{2pt}p$ higher hypergeometric series\medskip}

In this section new higher hypergeometric series for $\log\,p$ are obtained by
applying $W\hspace*{-2pt}Z_{s,t}$ transformations on Eqs.(\ref{6}--\ref{7}).
This convention is used where sums start from $n=0$
\begin{equation}
\log\hspace*{1pt}p=\frac{1}{c(p)}\cdot%
{\displaystyle\sum\limits_{n=0}^{\infty}}
\frac{P(n,p)}{R(n)}\cdot\mathcal{\rho(}p)^{n}\mathcal{\cdot}\left[
\begin{array}
[c]{cccc}%
r_{1} & r_{2} & ... & r_{d}\\
q_{1} & q_{2} & ... & q_{d}%
\end{array}
\right]  _{n}\tag{63}\label{63}%
\end{equation}
convergent for some region $\mathcal{D\subseteq%
\mathbb{C}
}$.\vspace*{-1pt} $c(p)$ is a simple rational function, $P(n,p)$ a polynomial
in $n$ of degree $d-1$ with polynomial coefficients and $\rho(p)$ is the
rational convergence rate function satisfying $|\mathcal{\rho(}p)|<1$ for
$p\in\mathcal{D}.$ Among many others, two outstanding cases are,$\vspace
*{0pt}$

\subsection{$\boldsymbol{d=4}$}

$W\hspace*{-2pt}Z_{2,3}^{\text{{\tiny (1)}}}|_{p}$ or $W\hspace*{-2pt}%
Z_{3,2}^{\text{{\tiny (2)}}}|_{p}$ transformations give this \vspace
*{-1pt}proven series convergent in $\mathcal{D=\{}p\,|\,p\in%
\mathbb{R}
_{>0}\wedge p_{0}^{-1}<p<p_{0}=$ {\small 28.52151255...}$\}$\ with rate
$\mathcal{O}((p-1)^{10})$ as $p\rightarrow1$,\vspace*{-2pt}%

\begin{equation}
\log\hspace*{1pt}p=-\frac{p-1}{p^{2}(p+1)^{5}}\cdot%
{\displaystyle\sum\limits_{n=0}^{\infty}}
\frac{P(n,p)}{R(n)}\left(  \frac{3^{3}}{2^{2}\,5^{5}}\cdot\dfrac{\,(p-1)^{10}%
}{\,p^{2}(p+1)^{6}}\right)  ^{n}\left[
\begin{array}
[c]{cccc}%
\genfrac{}{}{0pt}{0}{{}}{{}}%
1 & \frac{1}{2} & \frac{1}{6} & \frac{5}{6}%
\genfrac{}{}{0pt}{0}{{}}{{}}%
\\%
\genfrac{}{}{0pt}{0}{{}}{{}}%
\frac{1}{10} & \frac{3}{10} & \frac{7}{10} & \frac{9}{10}%
\genfrac{}{}{0pt}{0}{{}}{{}}%
\end{array}
\right]  _{n}\vspace*{4pt}\tag{64}\label{64}%
\end{equation}%
\[%
\begin{tabular}
[c]{ll}%
{\small \textit{P}(\textit{n,p}) \thinspace=} & {\small \ 216\thinspace
(\textit{p}}$^{2}${\small \thinspace+\thinspace8\textit{p}\thinspace
+\thinspace1)(27\textit{p}}$^{6}${\small \thinspace--\thinspace702\textit{p}%
}$^{5}${\small \thinspace--\thinspace1835\textit{p}}$^{4}${\small \thinspace
--\thinspace2980\textit{p}}$^{3}${\small \thinspace--\thinspace1835\textit{p}%
}$^{2}${\small \thinspace--\thinspace702\textit{p}\thinspace+\thinspace
27\thinspace)\thinspace\textit{n}}$^{3}\ ${\small +}$%
\genfrac{}{}{0pt}{0}{\genfrac{}{}{0pt}{1}{{}}{{}}}{{}}%
$\\
& {\small +\thinspace4\thinspace(81\textit{p}}$^{8}${\small \thinspace
--1674\textit{p}}$^{7}${\small \thinspace--27536\textit{p}}$^{6}%
${\small \thinspace--70486\textit{p}}$^{5}${\small \thinspace
--104770\textit{p}}$^{4}${\small \thinspace--70486\textit{p}}$^{3}%
${\small \thinspace--27536\textit{p}}$^{2}${\small \thinspace--1674\textit{p}%
\thinspace+\thinspace81\thinspace)\thinspace\textit{n}}$^{2}%
\genfrac{}{}{0pt}{0}{\genfrac{}{}{0pt}{1}{{}}{{}}}{{}}%
$\\
& {\small +\thinspace2\thinspace(69\textit{p}}$^{8}${\small \thinspace
--1674\textit{p}}$^{7}${\small \thinspace--30752\textit{p}}$^{6}%
${\small \thinspace--83670\textit{p}}$^{5}${\small \thinspace
--123146\textit{p}}$^{4}${\small \thinspace--83670\textit{p}}$^{3}%
${\small \thinspace--30752\textit{p}}$^{2}${\small \thinspace--1674\textit{p}%
\thinspace+\thinspace69\thinspace)\thinspace\textit{n}}$%
\genfrac{}{}{0pt}{0}{\genfrac{}{}{0pt}{1}{{}}{{}}}{{}}%
$\\
& {\small +\thinspace3\thinspace(5\textit{p}}$^{8}${\small \thinspace
--130\textit{p}}$^{7}${\small \thinspace--3184\textit{p}}$^{6}$%
{\small \thinspace--9694\textit{p}}$^{5}${\small \thinspace--14314\textit{p}%
}$^{4}${\small \thinspace--9694\textit{p}}$^{3}${\small \thinspace
--3184\textit{p}}$^{2}${\small \thinspace--130\textit{p}\thinspace
+\thinspace5\thinspace)}$\vspace*{6pt}%
\genfrac{}{}{0pt}{0}{\genfrac{}{}{0pt}{1}{{}}{{}}}{{}}%
$\\
{\small \textit{R}(\textit{n}) \ \ \thinspace=} & {\small \ \thinspace
20\thinspace(10\textit{n\thinspace}+\textit{\thinspace}%
1)(10\textit{n\thinspace}+\textit{\thinspace}3)(10\textit{n\thinspace
}+\textit{\thinspace}7)(10\textit{n\thinspace}+\textit{\thinspace}%
9)\vspace*{-4pt}}$%
\genfrac{}{}{0pt}{0}{\genfrac{}{}{0pt}{1}{{}}{{}}}{{}}%
$%
\end{tabular}
\]
$\vspace*{-4pt}$

\subsection{$\boldsymbol{d=6}$}

$W\hspace*{-2pt}Z_{4,3}^{\text{{\tiny (1)}}}|_{p}$ or $W\hspace*{-2pt}%
Z_{3,4}^{\text{{\tiny (2)}}}|_{p}$ transformations give this \vspace
*{-1pt}proven series convergent in $\mathcal{D=\{}p\,|\,p\in%
\mathbb{R}
_{>0}\wedge p_{0}^{-1}<p<p_{0}=$ {\small 34.808179...}$\}\ $with rate
$\mathcal{O}((p-1)^{14})$ as $p\rightarrow1$, \vspace*{0pt}where $P(n,p)$
is palindromic in $p$, so that the missing polynomial lower terms %
below must be completed by symmetry,\vspace*{4pt}%
\begin{equation}
\log\hspace*{1pt}p=\frac{1}{c(p)}\cdot%
{\displaystyle\sum\limits_{n=0}^{\infty}}
\frac{P(n,p)}{R(n)}\left(  \frac{3^{3}}{7^{7}}\cdot\dfrac{(p-1)^{14}}%
{\,p^{4}(p+1)^{6}}\right)  ^{n}\left[
\begin{array}
[c]{cccccc}%
\genfrac{}{}{0pt}{0}{{}}{{}}%
1 & \frac{1}{2} & \frac{1}{4} & \frac{3}{4} & \frac{1}{6} & \frac{5}{6}%
\genfrac{}{}{0pt}{0}{{}}{{}}%
\\%
\genfrac{}{}{0pt}{0}{{}}{{}}%
\frac{1}{14} & \frac{3}{14} & \frac{5}{14} & \frac{9}{14} & \frac{11}{14} &
\frac{13}{14}%
\genfrac{}{}{0pt}{0}{{}}{{}}%
\end{array}
\right]  _{n}\tag{65}\label{65}%
\end{equation}%
\[%
\begin{tabular}
[c]{ll}%
{\small \textit{c}(\textit{p}) \ \ \thinspace\ \thinspace=} &
{\small --\thinspace\textit{p}}$^{4}${\small (\textit{p}\thinspace
+\thinspace1)}$^{5}${\small /(\textit{p}\thinspace--\thinspace1)}$%
\genfrac{}{}{0pt}{0}{\genfrac{}{}{0pt}{1}{{}}{{}}}{\genfrac{}{}{0pt}{0}{{}%
}{{}}}%
$\\
{\small \textit{P}(\textit{n,p}) \thinspace=} & {\small 128\thinspace
(\textit{p}}$^{{\small 2}}${\small \thinspace+\thinspace5\textit{p}%
\thinspace+\thinspace1)(27\textit{p}}$^{{\small 10}}${\small \thinspace
--\thinspace648\textit{p}}$^{{\small 9}}${\small \thinspace+\thinspace
8208\textit{p}}$^{{\small 8}}${\small \thinspace--\thinspace74682\textit{p}%
}$^{{\small 7}}${\small \thinspace--\thinspace264859\textit{p}}$^{{\small 6}}%
${\small \thinspace--\thinspace411740\textit{p}}$^{{\small 5}}${\small ...
)\thinspace\textit{n}}$^{{\small 5}}%
\genfrac{}{}{0pt}{0}{\genfrac{}{}{0pt}{1}{{}}{{}}}{{}}%
$\\
\multicolumn{2}{l}{{\small +\thinspace32\thinspace(\thinspace270\textit{p}%
}$^{12}${\small \thinspace--5265\textit{p}}$^{11}${\small \thinspace+
53757\textit{p}}$^{10}${\small \thinspace--400383\textit{p}}$^{9}%
${\small \thinspace--7320270\textit{p}}$^{8}${\small \thinspace
--21251920\textit{p}}$^{7}${\small \thinspace--30355514\textit{p}}$^{6}%
${\small ... )\thinspace\textit{n}}$^{4}%
\genfrac{}{}{0pt}{0}{\genfrac{}{}{0pt}{1}{{}}{{}}}{{}}%
$}\\
\multicolumn{2}{l}{{\small +\thinspace8\thinspace(\thinspace1005\textit{p}%
}$^{12}${\small \thinspace--20040\textit{p}}$^{11}${\small \thinspace
+\thinspace215166\textit{p}}$^{10}${\small \thinspace--1773744\textit{p}}%
$^{9}${\small \thinspace--32182333\textit{p}}$^{8}${\small \thinspace
--94707848\textit{p}}$^{7}${\small \thinspace--135214940\textit{p}}$^{6}%
${\small ... )\thinspace\textit{n}}$^{3}%
\genfrac{}{}{0pt}{0}{\genfrac{}{}{0pt}{1}{{}}{{}}}{{}}%
$}\\
\multicolumn{2}{l}{{\small +\thinspace4\thinspace(\thinspace855\textit{p}%
}$^{12}${\small \thinspace--17340\textit{p}}$^{11}${\small \thinspace
+\thinspace195534\textit{p}}$^{10}${\small \thinspace--1831956\textit{p}}%
$^{9}${\small \thinspace--33025527\textit{p}}$^{8}${\small \thinspace
--98963632\textit{p}}$^{7}${\small \thinspace--141374300\textit{p}}$^{6}%
${\small ... )\thinspace\textit{n}}$^{2}%
\genfrac{}{}{0pt}{0}{\genfrac{}{}{0pt}{1}{{}}{{}}}{{}}%
$}\\
\multicolumn{2}{l}{{\small +\thinspace2\thinspace(\thinspace327\textit{p}%
}$^{12}${\small \thinspace--6708\textit{p}}$^{11}${\small \thinspace
+\thinspace78486\textit{p}}$^{10}${\small \thinspace--858504\textit{p}}$^{9}%
${\small \thinspace--15564079\textit{p}}$^{8}${\small \thinspace
--47782020\textit{p}}$^{7}${\small \thinspace--68437468\textit{p}}$^{6}%
${\small ... )\thinspace\textit{n}}$%
\genfrac{}{}{0pt}{0}{\genfrac{}{}{0pt}{1}{{}}{{}}}{{}}%
$}\\
\multicolumn{2}{l}{{\small +\thinspace3\thinspace(\thinspace15\textit{p}%
}$^{12}${\small \thinspace--310\textit{p}}$^{11}${\small \thinspace
+\thinspace3720\textit{p}}$^{10}${\small \thinspace--46330\textit{p}}$^{9}%
${\small \thinspace--887599\textit{p}}$^{8}${\small \thinspace
--2811664\textit{p}}$^{7}${\small \thinspace--4047184\textit{p}}$^{6}%
${\small ... )\thinspace}$%
\genfrac{}{}{0pt}{0}{\genfrac{}{}{0pt}{1}{{}}{{}}}{{}}%
\vspace*{6pt}$}\\
{\small \textit{R}(\textit{n}) \ \ \thinspace=} & {\small \ \thinspace
56\thinspace(14\textit{n\thinspace}+\textit{\thinspace}%
1)(14\textit{n\thinspace}+\textit{\thinspace}3)(14\textit{n\thinspace
}+\textit{\thinspace}5)(14\textit{n\thinspace}+\textit{\thinspace
}9)(14\textit{n\thinspace}+\textit{\thinspace}11)(14\textit{n\thinspace
}+\textit{\thinspace}13)}$%
\genfrac{}{}{0pt}{0}{\genfrac{}{}{0pt}{1}{{}}{{}}}{{}}%
\vspace*{4pt}$%
\end{tabular}
\]
Eqs.(\ref{64}--\ref{65}) can be also proven by BI method with
({\small \textit{m,}}$\nu$){\small \ = }({\small 3,2}) and ({\small 3,4})
respectively. Notice that in Eq.(\ref{18z}), Eq.(\ref{18m}) and Eqs.(\ref{64}%
--\ref{65}), $\log\,p$ can be converted into arc$\tanh(\zeta)$ and
$\arctan(\zeta)$ by means of elementary identities giving new highly efficient
hypergeometric series for these functions as well, however this research line
will be explored elsewhere.

\section{PARI GP Modules\bigskip}

Since there is no option of accompanying ancillary files besides the source TeX file, 
PARI GP scripts \texttt{myPSLQ} and \texttt{pfbeta}, among others, have been embedded as 
a large comment towards the end of the parent TeX file. This file shall be downloaded 
to extract (copy-paste) the code. Directions to install and test the scripts are placed at the 
beginning of such comment.\vspace*{-8pt}

\section{Acknowledgements\bigskip}

My acknowledgments to Jesús Guillera for valuable talks and comments leading me to search and find these powerful hypergeometric formulas. To Henri Cohen and all the PARI GP team, there would have been very difficult to get these fast identities without this 
highly efficient parallelized computing platform. Special acknowledgements to Alexander Yee for his powerful hypergeometric numerical computing y-cruncher software. Special thanks to Brian Beeler and Jordan Ranous at StorageReview Labs for providing high performance 128-core facilities to search and find Eqs.(\ref{8}$-$\ref{8a}$-$\ref{8b}) and Eq.(\ref{13}) and also to extend the decimal places known of \texttt{log(2)} up to \texttt{3e12} digits using the discovered Eq.(\ref{8} series. Special thanks as well to Lorenz Milla, Joshua Swanson and Bruce Sagan for their endorsement to this work.

\vfill\eject

\end{document}